\documentclass{article}

\usepackage{arxiv}

\usepackage[utf8]{inputenc} 
\usepackage[T1]{fontenc}    
\usepackage{hyperref}       
\usepackage{url}            
\usepackage{booktabs}       
\usepackage{amsfonts}       
\usepackage{nicefrac}       
\usepackage{microtype}      
\usepackage{lipsum}
\usepackage{subcaption}
\usepackage[final]{pdfpages}
\usepackage{array}
\newcolumntype{P}[1]{>{\centering\arraybackslash}p{#1}}
\usepackage{mathtools}
\usepackage{amsmath}
\usepackage{bm}

\DeclareMathOperator*{\argmin}{arg\,min}
\newcommand\norm[1]{\left\lVert#1\right\rVert}
\usepackage[super]{nth}
\usepackage{amssymb}
\numberwithin{equation}{section}
\usepackage{bookmark}
\usepackage{verbatim}
\usepackage{graphicx}
\usepackage[numbers]{natbib}
\usepackage{mwe}
\usepackage{color}

\newcommand{\?}{\:\!}
\newcommand{\T}{^{\?\mathrm{\scriptsize{T}}}}

\title{A Local Basis Approximation Approach for Nonlinear Parametric Model Order Reduction}

\author{Konstantinos Vlachas \\
    Dept. of Civil, Environmental and Geomatic Engr.\\
    ETH Zurich\\
    Zurich, Switzerland\\
    \texttt{vlachas@ibk.baug.ethz.ch} 
    \And
    Konstantinos Tatsis \\
    Dept. of Civil, Environmental and Geomatic Engr.\\
    ETH Zurich\\
    Zurich, Switzerland\\
    \texttt{tatsis@ibk.baug.ethz.ch}\\
    \And
    Konstantinos Agathos \\
    Dept. of Civil, Environmental and Geomatic Engr.\\
    ETH Zurich\\
    Zurich, Switzerland\\
    \texttt{agathos@ibk.baug.ethz.ch}
    \And
    Adam R. Brink \\
    Solid Mechanics\\
    Sandia National Laboratories\\
    Albuquerque, New Mexico,\\
    \texttt{arbrink@sandia.gov}
    \And
    Eleni Chatzi \\
    Dept. of Civil, Environmental and Geomatic Engr.\\
    ETH Zurich\\
    Zurich, Switzerland\\
    \texttt{chatzi@ibk.baug.ethz.ch}
}

\begin{document}
\maketitle

\begin{abstract}
The efficient condition assessment of engineered systems requires the coupling of high fidelity models with data extracted from the state of the system `as-is'. In enabling this task, this paper implements a parametric Model Order Reduction (pMOR) scheme for nonlinear structural dynamics, and the particular case of material nonlinearity. A physics-based parametric representation is developed, incorporating dependencies on system properties and/or excitation characteristics. The pMOR formulation relies on use of a Proper Orthogonal Decomposition applied to a series of snapshots of the nonlinear dynamic response. A new approach to manifold interpolation is proposed, with interpolation taking place on the reduced coefficient matrix mapping local bases to a global one. We demonstrate the performance of this approach firstly on the simple example of a shear-frame structure, and secondly on the more complex 3D numerical case study of an earthquake-excited wind turbine tower. Parametric dependence pertains to structural properties, as well as the temporal and spectral characteristics of the applied excitation. The developed parametric Reduced Order Model (pROM) can be exploited for a number of tasks including monitoring and diagnostics, control of vibrating structures, and residual life estimation of critical components.
\end{abstract}

\keywords{As-Built As-Deployed Structures \and  Parametric Model Order Reduction (pMOR) \and  Nonlinear Reduction \and  Reduced Bases Interpolation}

\section{Introduction} \label{Intro}

The emergence of digital twins as a main enabler for virtualization mandates the coupling of high fidelity simulations with data extracted from monitored systems \citep{citeE}. A `twin' of an operating system aims at precise representation of its response across the range of the system's regular and extreme loading conditions. This allows for robust design and effective diagnostics under uncertainty \citep{citeE3,citeE2}. However, the requirement for increased simulation accuracy implies a higher demand on computational resources, often compromising efficiency. For example, models that can capture local phenomena require the inclusion of higher-end detail during the modeling stage, which leads to complex numerical representations. To alleviate the burden of computation, a trade-off between accuracy and efficiency - tailored to the needs of the implementation - is needed. A means for maintaining accuracy lies in imprinting the underlying physics into the numerical representation. Efficiency is tackled by deriving low-dimensional models, which are capable of rapid computation, while sufficiently approximating the underlying high fidelity representation \citep{Benner,bookW}. This notion is referred to as Model Order Reduction and the derived low-dimensional model as Reduced Order Model (MOR and ROM respectively). For the sake of simplicity the ROM acronym is used throughout this paper to indicate the Reduced Order Model and MOR for the process of Model Order Reduction respectively. 

The available MOR literature spans from works in the domains of fluid dynamics \citep{Ballarin} and biomedical engineering \citep{Niroomandi}, to the fields of fracture mechanics \citep{Kerfriden2}, structural dynamics \citep{Amsallem}, monitoring and state estimation \citep{tat} in a more civil engineering oriented perspective. Depending on the domain of implementation and on the natural complexity of the addressed problem (linear, nonlinear or chaotic), several approaches and methodologies have been proposed. A comprehensive overview, serving as an initial basis, may be found in the works of \citet{bookW} and \citet{Antoulas}. \citet{CompareMOR} further offer a review and cross-comparison of established MOR techniques across scientific disciplines including structural dynamics, control and applied mathematics.

Under deterministically prescribed loads and system properties, a single evaluation of the ROM at an example sample represents a specific instance of the physical system under study. In this case, the ROM may only reliably describe the system's response for a given configuration, and around a narrow range of system parameters. In generalizing the applicability of the ROM formulation, it is necessary to express the delivered representation in terms of the parameters that enter the governing equations. This generalization is achieved via adoption of parametric Reduced Order Models (pROMs) \citep{Benner}. A dynamic pROM, in particular, should allow adequate approximation of the dynamics of the high fidelity system throughout the range of interest of modeling parameters. 

Assuming availability of a High Fidelity Model (HFM) of the system under study, the derivation of a ROM with parametric dependency usually requires querying this representation over multiple parametric inputs. An extensive overview of this class of intrusive methods can be found in \citet{Benner,bookW}. In the same context, data-driven and indirect, non intrusive methods have also been proposed. Although this extends beyond the focus of this paper, the interested reader is encouraged to refer to \citet{Mignolet} for a review on non intrusive methods and to \citet{add3} for suitable techniques to address nonlinearity. Regarding data-driven methods, \citep{Data2} can be treated as a starting point review. 

Although several approaches for pROMs do exist, not all of them are suitable for reproducing the time-domain response and the underlying physics of parametric, high order, nonlinear dynamical systems, represented by Finite Element (FE) models. An overview of relevant literature \citep{bookW,Benner,add2}, reveals the Proper Orthogonal Decomposition (POD - \citep{pod}) as the dominant reduction method for this class of problems. A short introductory review on the POD method for MOR can be found in \citet{Chinesta}, whereas the use of POD for identifying structures with geometric nonlinearities is reviewed in \citet{add1}. The POD technique relies on decomposition of time-series response data from training configurations of the system in order to identify a basis for a subspace of lower dimensions, where the response lies and where solutions are sought during the validation phase of the method. Since, typically, this subspace is of a much lower dimension than the HFM, this allows to greatly reduce the size of the problem. In principle this response could be obtained from either simulated or actually measured data. In the latter case a direct link is offered on fusion of structural models with monitoring data extracted from operating structures.

The POD approach is normally coupled with projection-based reduction methodologies \citep{Chinesta}. For example, a global reduction basis may be assembled for the entire parametric domain, which seems to be the method of choice when linear problems are addressed \citep{Amsallem}. In \citet{Agathos2019f} cracked, two dimensional solids are parameterized with respect to the geometric properties of the crack. The pROM is shown to efficiently approximate the static response for cracks that do not lie within the training set. In the same context, \citet{Creixell} propose an adaptive technique, where the pROM is embedded within an optimization process regarding the number of full order model evaluations. Similar considerations on reducing the required number of the training samples are discussed in \citet{Sample}. The global approach may be of further use when the phenomena under consideration are of localized nature. In the domain of fracture mechanics for instance, \citet{Kerfriden} proposed a ROM based on domain partitioning, focusing the numerical effort on the local domain of the defect. The structure domain is partitioned and reduced except for the localized domain of nonlinearity where the full model is evaluated.

The global approach might however yield inaccurate results or become computationally inefficient for nonlinear systems, where the response is strongly dependent on the input characteristics \citep{AmsalHP,Zimmer}. This implies that each new training parametric configuration contributes substantially new information to the POD basis. To address this issue, \citet{Amsal} assembled a pool of local POD bases instead of using a single global projection basis. The notion of locality refers to "neither space nor time a priori, but to the region of the manifold where the solution lies at a given parametric input or time instance" \citep{Amsallem}. In \citet{Amsal3} local subspaces are assembled with respect to the parametric domain of input, whereas in \citet{AmsalHP} with respect to time. 

To reproduce the response of the system in an `unseen' configuration, interpolation techniques are employed. In \citet{Amsal} the local projection bases are interpolated directly, with interpolation performed in the space tangent to the Grassmann manifold. This is done to ensure the interpolated bases maintain orthogonality properties. A recent study discussing detailed aspects of this approach is given in \citet{Zimmer2}. An alternative approach has been implemented in \citet{Degroot,Panzer2,linearp}, where interpolation of the projected ROM matrices is used. This is referred to as matrix interpolation and linear parametric ROMs are interpolated either in the original domain or in the space tangent to the Grassmann manifold \citep{Amsal}, after some form of coordinate transformation has been applied for consistency purposes \cite{linearp}. In this context, the studied interpolation schemes include spline interpolation \citep{Degroot} or Lagrange polynomials \citep{linearp}.

Several contributions across various scientific disciplines are based on the former two local interpolation notions, namely matrix and local bases interpolation. In \citet{Amsallem} for example, local bases are combined with k-means clustering-based interpolation algorithms. A similar approach is used in \citet{AmsalHP} where locality is defined with respect to time and a clustering approach with a suitable error metric for online updating is implemented. In \citet{loadpaths} Gaussian Processes and Bayesian strategies are adopted to address the optimal sampling of parametric input vectors, whereas similar considerations on optimizing the number of training samples while using a local bases interpolation approach are made in \citet{Taine,Soll,sample1}. A novel pMOR approach, suitable for contact in multibody dynamics is suggested in \citet{GearContact}. Global contact shapes are derived and the pROM is augmented to increase accuracy. Multi-parametric approximations are discussed in \citet{Washabaugh}, where ROMs are used to estimate steady-state flows over complex parameterized geometries. Another important contribution is \citet{sparsePOD} where an improved selection strategy for the basis vectors spanning the solution subspace is proposed. In \citet{CompareMOR2} a comparison of interpolation-based reduction techniques in the context of material removal in elastic multibody systems is discussed. Further contributions span across the domains of shape and design optimization \citep{Design}, solid and contact mechanics \citep{ContactBala,multilevel}, constrained optimization problems \citep{ZahrPDE}, highly nonlinear fluid-structure interaction problems \citep{Bala2} or coupling MOR with Component Mode Synthesis \citep{Recent2}.

The mentioned MOR methods focus mainly on identifying bases, able to represent the HFM solutions over a range of parameters. However, a main limiting factor for the performance of nonlinear projection-based ROMs is the actual projection of the system vectors and matrices containing nonlinear terms in the reduced space \citep{ECSW,aDEIM}. Therefore, a necessary component for formulating computationally efficient ROMs for FE-based formulations lies in the use of techniques to reduce the complexity of such projections, collectively referred to as hyper-reduction techniques. These second-level approximation strategies typically rely on evaluating the necessary projections on a limited, appropriately selected, set of elements in the FE mesh. Methods of this category include the Gauss–Newton with Approximated Tensors (GNAT) method presented in \citet{GNAT1}, the Discrete Empirical Interpolation Method (DEIM) suggested in \citet{DEIM} and the Energy-Conserving Sampling and Weighting hyper reduction method proposed in \citet{ECSW}. The computational efficiency of these approaches has been successfully demonstrated in \citet{GNAT2,tiso,Negri,AmsalHP,CompareMOR3}. Based on the considerations made in \citet{HPcompare} regarding numerical stability, structure preserving properties and overall efficiency, this paper implements the ECSW approach.

This paper focuses on development and demonstration of a pROM scheme for approximation of the time history response of structural systems featuring material nonlinearity, tied to phenomena such as plasticity and hysteresis. The topic of material nonlinearity has been addressed in the context of thermal loads for example in \citet{stressROM}, or with respect to impact analysis \citep{Taine,ECSW}. Our paper adds to this literature by implementing a pROM, able to model material nonlinearity and address parametric dependency pertaining in either the system or the excitation configuration. Specifically, we demonstrate the relevant aspects of parameterization in terms of influencing structural traits, such as material properties and hysteresis parameters, as well as in terms of parameterization of acting loads, e.g. earthquakes, in the temporal and spectral sense. Moreover, we implement a variant technique to the established local bases interpolation method in the sense that the dependence of a HFM on the characteristics of the configuration and/or loading parameters is expressed on a separate level to that of the snapshot procedure. The resulting pROM allows for accelerated computation, which is particularly critical for applications in monitoring and diagnostics, control of vibrating structures, and residual life estimation of critical components.

\section{Problem Statement and pROM Formulation} \label{PrSt}

Structural systems and components are exposed to harsh operating environments and adverse loadings. The resulting dynamics can be of nonlinear nature, albeit typically simplistically approximated by linear domain assumptions, which are oftentimes insufficient. A further complexity lies in the fact that the dynamics of real-life structural systems is dependent on material, geometrical, environmental or loading parameters. Therefore, when desired to reliably approximate the dynamics of an `as-built as-deployed' system, i.e., of a system in its operating state the aspects of nonlinearity and parametric variability need to be efficiently tackled. Naturally, this poses a limiting factor from the computational point of view, since the precise approximation of a system across the whole range of its defining parameters would require high fidelity simulations across multiple realizations of its defining parameter vector. This implies a high computational toll that prohibits adoption into practice, particularly for scenarios which necessitate a fast reaction, such as fault diagnosis and control. The aim of this work is to address the accurate modeling of `as built as-deployed' nonlinear dynamic structural systems (or components). The task is to deliver a sufficient reduced order approximation under variability of the structural properties or of the acting loads.

This section offers a short overview of the background knowledge required for the formulation of the proposed pROM. Firstly, the aspect of parametric dependency is visited by providing an overview of existing pROM formulations. Then, hyper-reduction is discussed as a further necessary reduction method for reducing the prohibitive computational toll associated with large degree of freedom finite element representations. 

\subsection{POD-based pROMs for nonlinear systems}\label{cur}
This subsection presents the fundamental background required for the treatment  of parametric dependencies within the MOR framework, which serves as the basis for establishing a parametric reduced order model (pROM) formulation for nonlinear systems. Within this context, a nonlinear dynamical system whose configuration depends on $l$ parameters, contained in the parameter vector $\mathbf{p}=[p_1,...,p_{\mathrm{l}}]^{\mathrm{T}} \in \Omega \subset \mathbb{R}^l$ is considered. The vibration response of such system is described by the governing equations of motion:
\begin{equation}
    \mathbf{M}(\mathbf{p})\ddot{\mathbf{u}}(t) +   \mathbf{g}\left(\mathbf{u}(t), \dot{\mathbf{u}}(t), \mathbf{p}\right) = \mathbf{f}(t,\mathbf{p}),
    \label{general}
\end{equation}
\noindent
where $\mathbf{u}(t) \in \mathbb{R}^n$ represents the system displacement, $\mathbf{M}(\mathbf{p}) \in \mathbb{R}^{n \times n}$ denotes the mass matrix,  $\mathbf{f}(t, \mathbf{p}) \in \mathbb{R}^{n}$ represents the vector of externally applied loads and $n$ is the order of the system, which physically represents the number of degrees of freedom. The nonlinearity of the system lies in the restoring force term $\mathbf{g}\left(\mathbf{u}(t), \dot{\mathbf{u}}(t)\right) \in\mathbb{R}^{n}$, which represents the resisting or internal forces of the system due to internal stresses and strains and is further dependent, along with the mass matrix and the externally applied excitation, on the parameter vector $\mathbf{p}$. This dependency may represent different system configurations depending on the target application, such as damage scenarios, which may be reflected on the stiffness and/or damping and mass matrices, or varying boundary conditions, which are dictated by the excitation vector.

The goal of parametric ROMs is to generate an equivalent system of dimension $r$, such that $r << n$ and the underlying physics along with the parametric dependencies of interest are further retained. Given the dependence of the governing system equations, represented by Equation \eqref{general}, on the parameters $\mathbf{p}$, the reduction step is herein performed for a number of sample points $\mathbf{p}_{\mathrm{j}}$ for $j=1,2,\dots,N$ in the parameter space using a projection-based strategy. As such, the solution of Equation \eqref{general} for a certain parameter sample $\mathbf{p}_{\mathrm{j}}$ is attracted to a lower dimensional subspace $S \subset \mathbb{R}^{n}$, spanned by the set of orthonormal basis vectors $\mathbf{V}(\mathbf{p}_{\mathrm{j}}) = \left[\mathbf{v}_1(\mathbf{p}_{\mathrm{j}}),\mathbf{v}_2(\mathbf{p}_{\mathrm{j}}),...,\mathbf{v}_{\mathrm{r}}(\mathbf{p}_{\mathrm{j}})\right]$, according to
\begin{equation}
    \mathbf{u}(t) = \mathbf{V}(\mathbf{p}_{\mathrm{j}}) \mathbf{u}_{\mathrm{r}}(t),
    \label{modal}
\end{equation}
\noindent
where $\mathbf{V}(\mathbf{p}_{\mathrm{j}})\in\mathbb{R}^{n\times r}$ is called the projection basis and the reduced size vector $\mathbf{u}_{\mathrm{r}} \in \mathbb{R}^{r}$ defines the components of the solution in this basis. Thereafter, the reduced-order representation of the HFM is obtained by means of a Galerkin projection, which is carried out by substituting Equation \eqref{modal} into Equation \eqref{general}, pre-multiplying with the transpose of the projection basis $\mathbf{V}(\mathbf{p}_{\mathrm{j}})\T$ and lastly imposing the orthogonality condition of the occurring residuals with respect to the projection basis. This results in the following system of equations:
\begin{equation}
    \mathbf{M}_{\mathrm{r}}(\mathbf{p}_{\mathrm{j}})\ddot{\mathbf{u}}_{\mathrm{r}}(t) +   \mathbf{g}_{\mathrm{r}}\left(\mathbf{u}(t), \dot{\mathbf{u}}(t), \mathbf{p}_{\mathrm{j}}\right) = \mathbf{f}_{\mathrm{r}}(t,\mathbf{p}_{\mathrm{j}})
    \label{reducedlmor}
\end{equation}

\noindent
where $\mathbf{M}_{\mathrm{r}}(\mathbf{p}_{\mathrm{j}})\in \mathbb{R}^{r \times r}$ and $\mathbf{g}_{\mathrm{r}}\left(\mathbf{u}(t), \dot{\mathbf{u}}(t), \mathbf{p}_{\mathrm{j}}\right)\in \mathbb{R}^{r}$ denote the reduced-order mass matrix and restoring force vector respectively, while $\mathbf{f}_{\mathrm{r}}(t, \mathbf{p}_{\mathrm{j}})\in\mathbb{R}^{r}$ designates the generalized vector of external forces, as follows
\begin{align}
    \mathbf{M}_{\mathrm{r}}(\mathbf{p}_{\mathrm{j}}) & = \mathbf{V}(\mathbf{p}_{\mathrm{j}})\T\mathbf{M}(\mathbf{p}_{\mathrm{j}})\mathbf{V}(\mathbf{p}_{\mathrm{j}}) \nonumber \\[2mm]
    \mathbf{g}_{\mathrm{r}}(\mathbf{p}_{\mathrm{j}}) & = \mathbf{V}(\mathbf{p}_{\mathrm{j}})\T \mathbf{g}\left(\mathbf{u}(t), \dot{\mathbf{u}}(t), \mathbf{p}_{\mathrm{j}}\right) \nonumber \\[2mm]
    \mathbf{f}_{\mathrm{r}}(\mathbf{p}_{\mathrm{j}}) & = \mathbf{V}(\mathbf{p}_{\mathrm{j}})\T\mathbf{f}(t, \mathbf{p}_{\mathrm{j}})
    \label{matproj}
\end{align} 
\noindent
The computation of the orthonormal basis vectors $\mathbf{V}(\mathbf{p}_{\mathrm{j}}) = \left[\mathbf{v}_1(\mathbf{p}_{\mathrm{j}}),\mathbf{v}_2(\mathbf{p}_{\mathrm{j}}),...,\mathbf{v}_{\mathrm{r}}(\mathbf{p}_{\mathrm{j}})\right]$, which is the key element of the reduction step, is typically carried out by means of Proper Orthogonal Decomposition (POD) \citep{pod}. As such, a pool of displacement field samples is collected from the time history analysis of the HFM. Each one of the pools is extracted from a unique parameter configuration $\mathbf{p}_{\mathrm{j}}$ of the full order model, which is henceforth termed as snapshot. Thereafter, the information of all simulated pools of samples, or equivalently snapshots, is collected on a global basis $\mathbf{V}_{\mathrm{full}}$, which essentially captures the parametric dependence of the model, since response information across the entire parameter space is assembled through sampling. A Singular Value Decomposition (SVD) is lastly applied to $\mathbf{V}_{\mathrm{full}}$ in order to obtain the principal orthonormal components of the reduction basis $\mathbf{V}$ spanning the lower subspace $S$ of the solution. The error measure used in \citep{tiso} is herein employed for this purpose as well. These steps represent the offline phase of reduction, which produces the global SVD-based projection basis $\mathbf{V}$. The latter may be utilized for global system order reduction, according to Equation \eqref{matproj}, with the aim of reflecting the underlying dynamics and therefore approximating the response of the model at unseen parameter samples.

One of the limitations of such an approach lies on the fact that the response of nonlinear systems is strongly dependent on the parameter values and may be dominated by localized phenomena which are owed to the nonlinear terms. This implies that the response may well lie on spaces that cannot be spanned by a single global basis, unless the latter comprises a large number of basis vectors.

As a result, localized features that are observed only in a restricted parameter subdomain end up affecting and determining the overall estimation capabilities of the ROM \citep{Zimmer}. To address these issues, a different strategy can be employed, as was initially highlighted in section \ref{Intro}. Instead of formulating a global projection basis, a pool of local POD bases can be assembled during the offline phase, with each one of those bases corresponding to a unique sample, or family of parameter samples. To reproduce then the response of the system for an `unseen' configuration in an online manner, interpolation techniques can be employed with the aim of estimating the local POD basis at the required parameter point and subsequently obtain the reduced matrices based on Equation \eqref{matproj}.

Although computationally more efficient, the implementation of local basis interpolation schemes constitutes a rather non-trivial task. This is owed to the fact that each local basis is inherently characterized by specific algebraic properties that need be maintained upon interpolation. Specifically, the local bases are derived using POD on the response of the structural system and the projected matrices are then employed to integrate the system's response equations. Therefore, the orthogonality of the interpolated subspaces as well as the positive-definiteness of the occurring reduced-order matrices need be retained. As such, a trivial element-wise interpolation of the local bases would result in non-orthogonal basis vectors, leading to high approximation errors and stability issues, calling thus for more sophisticated interpolation schemes, able to account for the underlying constraints.

To alleviate this issue, interpolation is herein performed according to the approach described in the work of \citet{Amsal}. Within this context, each local basis is initially mapped to the tangent space of the Grassmann manifold, which essentially preserves the underlying orthogonality property. The mapped data is subsequently interpolated in order to yield the sought for basis at the unseen parameter point and the interpolated result is lastly mapped back to the manifold. The detailed mathematical formulation pertaining to the manifold mapping process and interpolation step can be found in \citep{Zimmer}. This projection space is selected based on the following considerations:
\begin{itemize}
    \item It is a well-established concept in the MOR community for constrained matrix and/or projection basis interpolation
    \item The mapping algorithms to and from the projection space are well-known and relatively straight-forward to implement
    \item The forward and backward manifold projections as well as the interpolation on this space, retains the  orthogonality property.
\end{itemize}

\subsection{Hyper-Reduction}\label{sechyper}
In the previous subsection, the parameter dependency was injected into the MOR framework with the aim of guaranteeing efficiency of the approximation across the entire parameter space, i.e., for different values of structural properties and different load characteristics. However, the aspect of computational efficiency is yet to be discussed since the computational complexity of the problem remains bounded by the size of the full order model. To this end, this subsection overviews the concept of hyper-reduction for further reducing the computational cost of the pROM.

During formulation of the ROM of dimension $r$ in Equation \eqref{reducedlmor}, it was specified that $r<<n$, i.e., the ROM constructed using a projection-based strategy features a significantly lower dimension than the original HFM. However, the same does not apply for the overall computational complexity of the problem. This is owed to the fact that the evaluation of the nonlinear contributions at the element-level for each time-step of the integration process dominates the complexity of the overall framework, with the associated computational toll scaling with the full-order dimension of the system $n$ \citep{AmsalHP,ECSW}. 

To address this bottleneck, the nonlinear pROM should be equipped with an additional level of reduction in order to remove the dependency of the computational procedure on the full-order model dimension $n$. This is accomplished with the so-called hyper-reduction strategy, which may be employed in a number of variants, which are already well-documented in the literature \citep{ECSW,GNAT1,GNAT2,DEIM,aDEIM}. Herein, the Energy-Conserving Sampling and Weighting (ECSW) scheme is implemented, which was initially introduced in \citet{ECSW}. This approach is chosen due to its physics-based formulation, which naturally suits for FE computations. Moreover, it features a stable performance and has been further shown to preserve the Lagrangian structure associated with Hamilton's principle \citep{HPcompare}. The fundamental aspects of the method are discussed below.

As a first step, a distinction is made between force vectors that can and cannot be efficiently evaluated online. No hyper-reduction is required for the former, while the latter can be written as a sum of elemental contributions, as follows
\begin{equation}
    \mathbf{g}_{\mathrm{r}}^{\bullet}\left(\mathbf{u}(t)\right) = \sum_{e=1}^{n_{e}} \mathbf{V}_\mathrm{e}^{\mathrm{T}}\mathbf{g}_{e}^{\bullet}\left(\mathbf{V}_\mathrm{e} \mathbf{u}(t)\right)
    \label{eq:reduced_forces}
\end{equation}
where, if $r$ denotes the dimension of the ROM, $\mathbf{V}_{\mathrm{e}} \in \mathbb{R}^{n_{\mathrm{e}} \times r}$ is the projection basis evaluated only at the degrees of freedom related with element $e$, $n_{\mathrm{e}}$ is the number of elements and $\mathbf{g}_{e}^{\bullet}$ denotes the element forces.

The forces of Equation \eqref{eq:reduced_forces} can be approximated via a weighted sum over a subset $\widetilde{E}$ of the total elements, so that
\begin{equation}
    \mathbf{g}_{\mathrm{r}}^{\bullet}\left(\mathbf{u}(t)\right) \approx \widetilde{\mathbf{g}}_{\mathrm{r}}^{\bullet}(\mathbf{u}(t)) = \sum_{\widetilde{E}} \xi^{\ast}_{\mathrm{e}} \mathbf{V}_{\mathrm{e}}^{\mathrm{T}} \mathbf{g}_{\mathrm{e}}^{\bullet}\left(\mathbf{V}_{\mathrm{e}} \mathbf{u}(t)\right)
    \label{ECSW}
\end{equation}
\noindent 
where $\xi^{\ast}_{\mathrm{e}}$ are non-negative weights assigned to the elements of subset $\widetilde{E}$. In this sense, the number of operations required to compute the reduced forces scales with the number $\widetilde{n}_{\mathrm{e}}$ of elements in $\widetilde{E}$ and not with the total number of elements. Therefore, the efficiency of the pROM can be drastically increased if the number of elements in $\widetilde{E}$ is such that $\widetilde{n}_{\mathrm{e}} \ll n_{\mathrm{e}}$.

Equation \eqref{eq:reduced_forces} indicates that if the components of each row of $\mathbf{V}_\mathrm{e}$ are treated as virtual displacements, then each entry of the vector $\mathbf{g}_{\mathrm{r}}^{\bullet}(\mathbf{u}(t))$ corresponds to the work produced by the forces $\mathbf{g}_{\mathrm{e}}^{\bullet}(\mathbf{V}_\mathrm{e} \mathbf{u}(t))$ over the whole mesh along these displacements. Moreover, this work can be written as the sum, over all elements, of the work produced in each element. Thus, the set $\widetilde{E}$ and weights $\xi^{\ast}_{\mathrm{e}}$ can be determined by minimizing the difference between the work produced by the reduced and full mesh over a set of training configurations $n_{\mathrm{s}}$. To this end, element contributions are stored in a matrix $\mathbf{G}$, while the total work over the entire meshed domain, for different configurations, is stored in the elements of vector $\mathbf{b}$:
\begin{gather}
G_{\mathrm{ie}} = \mathbf{V}_\mathrm{e}^T \mathbf{g}_{\mathrm{e}}^{\bullet}(\mathbf{V}_\mathrm{e} \mathbf{u}(t)) \nonumber \\
b_{\mathrm{i}} = \sum_{e=1}^{\mathrm{n_{\mathrm{e}}}} G_{\mathrm{ie}} \label{xi2}\\
\mathbf{G} =
\begin{bmatrix} 
G_{11} & \dots & G_{\mathrm{1n_{\mathrm{e}}}} \\
\vdots & \ddots & \\
G_{\mathrm{n_{\mathrm{s}}}1} & \dots & G_{\mathrm{n_{\mathrm{s}} n_{\mathrm{e}}}}
\end{bmatrix}
\quad
\mathbf{b} =
\begin{bmatrix} 
b_{1}\\
\vdots \\
b_{\mathrm{n_{\mathrm{s}}}}
\end{bmatrix}
\label{HPmatrixes}
\end{gather}
With these definitions in place, $\widetilde{E}$ and $\xi^{\ast}_{\mathrm{e}}$ can be obtained as the solution of the constrained minimization problem
\begin{gather}
    \bm{\xi}^{\ast} = \argmin_{\bm{\xi} \in \Phi} \norm{\bm{\xi}}_0\hspace{3mm}\textrm{with}\hspace{3mm}
    \Phi = \bm{\xi} \in \mathbb{R}^{n_{\mathrm{e}}}: \norm{\mathbf{G} \bm{\xi}-\mathbf{b}}_2 \leq \tau \norm{\mathbf{b}}_2,\,\, \xi_e \geq 0
    \label{xi}
\end{gather}
where $\norm{\bullet}_0$  and $\norm{\bullet}_2$ denote the $L_0$ and $L_2$ norms respectively, $\Phi$ represents the feasible set of candidate solutions, and $\tau$ denotes a user-defined tolerance regulating the accuracy of the approximation. The above problem can be efficiently solved using the sparse Non-Negative Least Squares (sparse NNLS) algorithm \citep{NLS,ECSW}.

\section{Local Basis Coefficients Interpolation for pMOR}\label{3}
Given the problem of Equation \eqref{general}, our goal is to construct an accurate and computationally efficient pROM that reproduces the underlying dynamic behavior of the system across a range of system configurations and/or loading parameters. In doing so, a local basis interpolation approach is developed herein. This approach is a new variant of the local bases interpolation technique introduced in \citep{Amsal}. To this end, a new element is introduced in the interpolation technique and the parametric dependence of the problem is expressed on a separate level to that of the snapshot procedure or the local subspaces. This section describes the details of the proposed approach and the overall algorithmic framework. 

As a first step, the parameter space of the problem is sampled and partitioned in a set of subdomains. These subdomains are defined based on the similarity of the response and the underlying dynamics of the samples. Samples with similar displacement time histories are grouped together and considered as a local region. In the numerical case studies of this paper, the principal independent components of the response time histories of neighboring parameter samples are assumed to span the same subspace, which further implies a high degree of similarity among the corresponding displacement time histories. In other words, similarity of snapshots implies a close distance metric on the parameter space. Herein, the domain is partitioned employing a simplified structured grid of rectangular regions, which is proven to  deliver satisfactory results. Once this grid is defined, the HFM is simulated for all the points in the grid, plus some additional points which are utilized for validation, with the corresponding pool of snapshots hence created.

The approach implemented in this paper attempts to combine and bridge, to a certain extent, the single global basis approach and the local bases interpolation strategy. Within this context, a pool of local bases is constructed, namely one subspace spanned by the principal components of the snapshot for each one of the training samples. These local bases are then interpolated to generate the projection basis for any unseen point of interest. As explained in subsection \ref{cur}, these local bases need to be projected to an appropriate space, i.e., the tangent space to the Grassmann manifold, capable of preserving the fundamental property of orthogonality for the interpolated quantities. An overview of these steps is given in Table \ref{table3}.

\begin{figure}[!ht]
	\centering
	\includegraphics[scale=0.70]{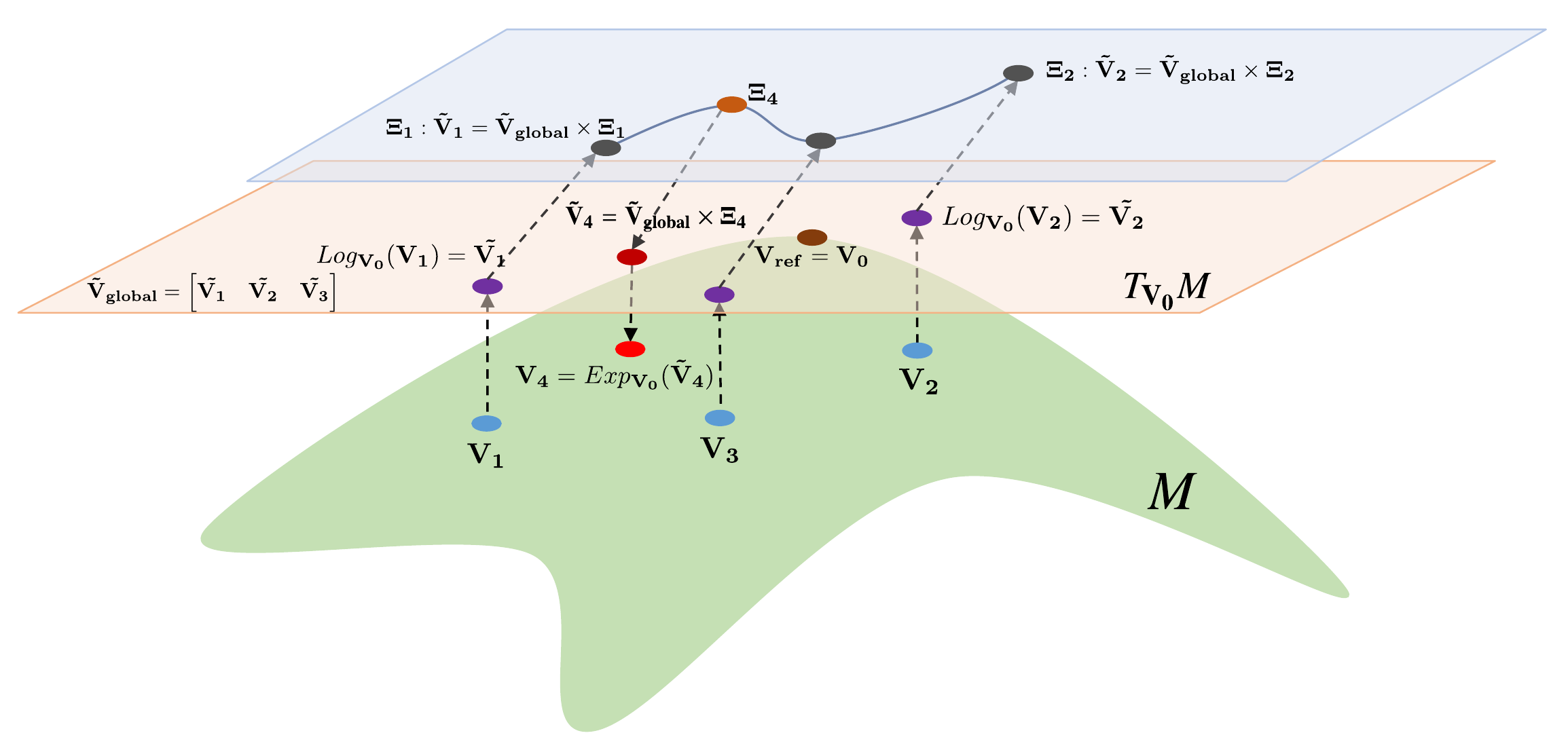}
	\caption{Graphical description of the pROM local bases interpolation approach implemented on this paper. The green surface represents the Grassmann Manifold ($M$) and the circular points denoted by $\mathbf{V}_{\mathrm{i}}, i=1,2,3$ the training samples. The validation parametric sample is denoted by $\mathbf{V}_4$ and the Exp and Log notations refer to the mapping processes to and from the manifold \citep{Zimmer}. The orange surface is the tangent space drawn to point $\mathbf{V_{\mathrm{ref}}}$. The representation is inspired from \citep{Amsal}.}
	\label{fig:manifold}
\end{figure}

To complete the definition of the tangent space projection strategy, a reference point needs to be defined. This represents the point where the tangent space to the Grassmann manifold is drawn. In Figure \ref{fig:manifold} this is represented by the brown point $\mathbf{V}_0$. The center point for every rectangular subdomain is considered herein as the reference point of each region for the definition of the tangent space. In this manner, we ensure that the subspaces of the remaining training samples of the subdomain are sufficiently close to the reference point, thus reducing projection errors due to curvature of the manifold \cite{Zimmer2}. After obtaining the pool of local bases and projecting them to a suitable tangent space, a `global' POD basis for every subdomain is assembled. As a result, we have constructed a single 'global' basis and a pool of local bases for every region of the parameter domain of interest. These first steps of the procedure are depicted graphically in Figure \ref{fig:manifold}. The green surface represents the Grassmann Manifold and the circular points denoted by $\mathbf{V}_{\mathrm{i}}, i=1,2,3$ the training samples of an example subdomain. The validation parameter sample is denoted by $\mathbf{V}_4$ and the Exp and Log operators refer to the mapping processes to and from the manifold \citep{Zimmer}.

After projecting the local bases to the tangent space, element-wise interpolation between the entries of the local bases could be adopted as a straightforward approach. This would imply interpolating between the points on the orange surface of Figure \ref{fig:manifold}, which would make the number of operation scale with the dimension $n$ of the HFM,  regardless of the interpolation scheme. Specifically, if each local basis $\mathbf{V}_{\mathrm{i}}$ retains $N_{\mathrm{m}}$ columns as principal modes, then $\mathbf{V}_{\mathrm{i}} \in \mathbb{R}^{n \times N_{\mathrm{m}}}$. The respective interpolation operations are $\mathcal{O}(n N_{\mathrm{m}})$. To get around this dependency, we adopt a different approach, which consists in projecting each local basis to the global POD basis of the region through a coefficient matrix $ \mathbf{\Xi} $. So, if $rg$ denotes the reference number of the region and $i$ the number of a snapshot in the region, the global basis is denoted by $\mathbf{V}^{\mathrm{rg}}_{\mathrm{global}}$ and the local bases as $\mathbf{V}^{\mathrm{rg}}_{\mathrm{local},i}$. The corresponding mathematical expression reads

\begin{equation}
    \mathbf{V}^{\mathrm{rg}}_{\mathrm{local},i} = \mathbf{V}^{\mathrm{rg}}_{\mathrm{global}} \mathbf{\Xi} _{\mathrm{i}}
    \label{linearcoeff}
\end{equation}

\noindent
with $\mathbf{V}^{\mathrm{rg}}_{\mathrm{local},i} \in \mathbb{R}^{n \times N_{\mathrm{m}}}$,  $\mathbf{V}^{\mathrm{rg}}_{\mathrm{global}} \in \mathbb{R}^{n\times (N_{\mathrm{s}} \times N_{\mathrm{m}})}$ and $ \mathbf{\Xi} \in \mathbb{R}^{(N_{\mathrm{s}} \times N_{\mathrm{m}}) \times N_{\mathrm{m}}}$, where $N_{\mathrm{s}}$ is the total number of generated snapshots and $N_{\mathrm{m}}$ is the number of independent components contained in each local basis, which are further assumed to be independent to the $N_{\mathrm{l}}$ components of the obtained local bases. In this manner, only an interpolation of the coefficient matrix per validation point $\mathbf{\Xi}$ is required in order to obtain a local projection basis. This interpolation scheme is graphically presented in Figure \ref{fig:manifold} by the light blue surface, which represents the domain spanned by the coefficient matrices, while the interpolation of the coefficient matrices is represented by the corresponding spline. Therefore, the proposed scheme requires interpolation of the coefficient matrix, which comprises only $\left(N_{\mathrm{s}} \times N_{\mathrm{m}}\right) \times N_{\mathrm{m}}$ entries, removing thus the dependency on the large dimension $n$ of the full problem. To the contrary, a straightforward interpolation scheme would need to estimate $n\times N_{\mathrm{m}}$ entries, which would be probably accomplished by $N_{\mathrm{m}}$ interpolation schemes, one for each $n$-length vector.

\begin{table}
	\caption{Reduction Framework Algorithmic Process. $N_{\mathrm{s}}$ denotes the number of training snapshots, $N_{dof}$ the degrees of freedom of the model,$N_{\mathrm{t}}$ the number of simulated timesteps and $N_{modes}$ the independent components of the reduction basis. Notation is kept similar to the Equations of sections \ref{PrSt} and \ref{3}}
	\noindent\rule[0.5ex]{\linewidth}{2pt}
	\label{table3}
	\begin{tabular} { p{13.5cm} }
	    \textbf{Notation:}\\
	     $N_{\mathrm{s}}$: Number of training samples, $N_{\mathrm{t}}$: Number of timesteps\\
	     $N_{\mathrm{dof}}$:Number of Degrees of Freedom of the model\\
	     $N_{\mathrm{modes}}$: Independent components/modes of a local reduction basis\\
	     $N_{\mathrm{k}}$: Number of samples per local region,
	     $N_{\mathrm{rg}}$: Number of local regions\\
	     $N_{\mathrm{reduced}}$: Independent components/modes of a global reduction basis\\
	     \hline
	     \\
		\textbf{Input:}\\ Parameter vector of training configurations $\mathbf{p}=\left[\mathbf{p}_1, \mathbf{p}_2, \dots ,\mathbf{p}_{N_{\mathrm{s}}}\right]$,\\ Validation configuration $\mathbf{p}_{\mathrm{q}}$ to be interpolated \\
		\textbf{Output:}\\
		Displacements of validation configuration estimated with ROM\\
		\hline
		\\
		\textbf{Of\-f\-line phase}\\
		1:\textbf{for} b=1,...,$N_{\mathrm{s}}$ \textbf{do}\\
		2:\quad Simulate full order model response for $\mathbf{p}=\mathbf{p}_{\mathrm{b}}$\\
		3:\quad Obtain displacement time history $\mathbf{U}_{\mathrm{b}} \in \mathbb{R}^{N_{\mathrm{dof}} \times N_{\mathrm{t}}}$\\
		4:\quad Perform SVD to obtain local basis $\mathbf{V}_{\mathrm{b}} \in \mathbb{R}^{N_{\mathrm{dof}} \times N_{\mathrm{modes}}}$\\
		5: \textbf{end for}\\
		6: Group* local regions based on samples' similarity of response and underlying dynamics\\
		7: \textbf{for} j=1,...,$N_{\mathrm{rg}}$ \textbf{do}:\\
		8: \quad Pick centroid of region ($j_0$) as reference to 'draw' the tangent plane $\mathbf{T}^{j}_{\mathbf{V}_{j_0}}$\\ 
		9: \quad Map matrices $\mathbf{V}^{\mathrm{j}}_{1,2,...,N_{\mathrm{k}}}$ to $ \widetilde{\mathbf{V}}^{\mathrm{j}}_{1,2,...,N_{\mathrm{k}}}$ spanning the tangent space $\mathbf{T}^{j}_{\mathbf{V}_{j_0}}$\\ \quad \quad \quad employing the Grassmann Matrix Logarithm \citep{Zimmer}\\
		10:\quad Assemble local bases $ \widetilde{\mathbf{V}}^{\mathrm{j}}_{1,2,...,N_{\mathrm{k}}}$ to a global matrix $\mathbf{V}^{\mathrm{j}}_{\mathrm{global}}$\\
		11:\quad Perform SVD to $\mathbf{V}^{\mathrm{j}}_{\mathrm{global}} \in \mathbb{R}^{N_{\mathrm{dof}} \times \left(N_{\mathrm{k}} \times N_{\mathrm{modes}} \right)}$ and obtain $\widetilde{\mathbf{V}}^{\mathrm{j}}_{\mathrm{global}} \in \mathbb{R}^{N_{\mathrm{dof}} \times N_{\mathrm{reduced}}}$\\
		12:\quad On the tangent space $\mathbf{T}^{j}_{\mathbf{V}_{j_0}}$, solve $\widetilde{\mathbf{V}}^{\mathrm{j}}_{1,2,...,N_{\mathrm{k}}} = \widetilde{\mathbf{V}}^{\mathrm{j}}_{\mathrm{global}} \times \mathbf{\Xi}_{1,2,...,N_{\mathrm{k}}}$ and store $\mathbf{\Xi}_{1,2,...,N_{\mathrm{k}}}$\\
		13:\textbf{end for}\\
		\\
		\textbf{Online phase}\\
		1*: Identify Region $rg$- Interpolate coefficients $\mathbf{\Xi}$ on the tangent space \\ \quad \quad to estimate $\mathbf{\Xi}_{\mathrm{q}}$ for parametric point $\mathbf{p}_{\mathrm{q}}$ \\
		2*: Compute local basis $\widetilde{\mathbf{V}}^{rg}_{\mathrm{q}} = \widetilde{\mathbf{V}}^{\mathrm{rg}}_{\mathrm{global}} \times \mathbf{\Xi}_{\mathrm{q}}$\\
		3*: Map local basis $\widetilde{\mathbf{V}}^{rg}_{\mathrm{q}}$ back to  ${\mathbf{V}}^{rg}_{\mathrm{q}}$ on the original parameter space employing the\\ \quad \quad Grassmann Matrix Exponent \citep{Zimmer}\\
		4*: Formulate the reduced order matrices based on Equation \eqref{matproj}\\
		5*: Simulate the reduced order model response (Equation \eqref{reducedlmor}) and obtain the \\ \quad \quad  displacement time history $\mathbf{U}_{\mathrm{q}}$\\
		\\
		*: Local regions are defined manually. The underlying assumption is that subspaces of close training snapshots can be interpolated to derive an accurate validation subspace for a new parametric point.   
	\end{tabular}
	\noindent\rule[0.5ex]{\linewidth}{2pt}
\end{table}

As briefly discussed in section \ref{PrSt}, the interpolation is carried out on the tangent space of the manifold and aims to exploit the fact that the solution in the vicinity of each parameter sample can be spanned by a local subspace, which is herein referred to as the local basis. Within this context, the local basis at an unseen point of the parameter space can be expressed as a combination of independent groups of basis vectors, with each one of those groups extracted from the solution obtained for a number of training points. In contrast with the element-wise interpolation among local bases, which is a rather numerical remedy, the considered interpolation scheme is capable of preserving in a sense the physical interpretability, provided that the solution at an unseen parameter point can be well spanned by a subspace similar to the ones spanning the solutions at the training points.

The detailed steps of the algorithm employed in this paper are documented in Table \ref{table3}, while a schematic representation of the projection and interpolation steps is provided in Figure \ref{fig:manifold}. All in all, in this paper a projection-based local bases interpolation scheme is formulated to express the dependence of a HFM on the characteristics of the system's configuration and/or loading parameters. We propose a new approach to this local bases interpolation, with the aforementioned dependence expressed on a separate level to that of the snapshot procedure, namely on the level of the global basis coefficients $\mathbf{\Xi}$. This approach is adopted in order to generate a pROM, able to accurately capture the underlying dynamics of the HFM across the parameter space. For validation purposes, the performance of the pROM is evaluated numerically in two cases studies featuring both nonlinearity and parametric dependence, with the latter pertaining to either the system properties or the characteristics of the external excitation. In addition, the proposed approach is compared with a global basis strategy and element-wise local bases interpolation on the tangent space. 

\section{Applications} \label{4}
After introducing the methodological background in section \ref{PrSt} and describing the proposed approach for pROM construction based on local bases interpolation in section \ref{3}, this section focuses on numerically verifying the efficiency of the proposed technique by means of two case studies. First, an academic example of a two story building with nonlinear links is simulated. The parametric dependency enters the configuration of the Bouc Wen hysteretic springs \citep{bw}, which are used to model the nonlinear couplings. This example demonstrates the potential of the pROM to handle nonlinear effects when the parametric dependency is introduced into the system material properties. Then, an actual scale wind turbine tower is simulated under earthquake excitation and under the assumption of material nonlinearity, with the parametric dependency in this case pertaining to the spectral characteristics of the excitation.

For the performance assessment of the pROMs, the parameter space is divided into rectangular subdomains, as discussed in section \ref{3}, and four different model configurations are examined. The first pROM uses a single global basis, which is extracted from training samples of the entire parameter domain while the second one utilizes a local basis technique for each subdomain of the parameter space. The other two pROMs employ local bases interpolation, with the first one performing element-wise interpolation on the entries of the projected to the tangent space of the Grassman manifold local bases, while the last one utilizes the proposed interpolation which was elaborated in section \ref{3}. The interpolation weights are computed based on the distance of the training configurations from the reference point of the corresponding parameter space subregion. The four different scenarios along with their corresponding description are summarized in Table \ref{tableref}. 

\begin{table}
	\caption{Reference table for compared pROMs.}
	\label{tableref}
	    \centering{
		\begin{tabular} { P{5cm} P{8cm}}
		\toprule
    		pROM Reference Name & Description\\
    		\hline
    		\\
    		Global Basis& Projection basis assembled using training snapshots across the whole parametric domain\\
    		Local Basis & Projection bases assembled using training snapshots of a subdomain\\
    		Entries Interpolation & Projection bases assembled by element-wise interpolation of local bases on the tangent space \\
    		Coefficients Interpolation & Projection bases assembled by coefficients interpolation on the tangent space based on formulation in section \ref{3}.\\
    	\bottomrule
	    \end{tabular}}
\end{table}
\noindent
To assess the accuracy of each pROM configuration, the following error norm is utilized, which constitutes a well established measure in the literature \cite{ECSW}

\begin{equation}
\mathbb{R}\mathbb{E}_{Q}=\frac{\sqrt{(Q_{\mathrm{HFM}}-{Q}_{\mathrm{ROM}})\T(Q_{\mathrm{HFM}}-{Q}_{\mathrm{ROM}})}}{\sqrt{(Q_{\mathrm{HFM}}\T{Q}_{\mathrm{ROM}})}},
\label{errornorm}
\end{equation}
where subscript $\mathrm{HFM}$ stands for High Fidelity Model results, $\mathrm{ROM}$ represents the ROM approximation and $Q$ denotes the quantity of interest. This can be either displacement, strains, stresses or restoring forces depending on the specific context of the implementation. Regarding displacements and restoring forces, Q denotes the time history response of the quantity of interest whereas for stresses this error measure is evaluated for the last time step of the analysis, namely the remaining values. 

In what concerns the computational cost, all simulations were performed in a single core on an Intel(R) Xeon(R) CPU machine and for the sake of fairness, the same integration time step is used for both pROM and HFM. The corresponding speed up factor is calculated as the ratio of CPU time required for an evaluation of the pROM over the time needed for the evaluation of the HFM, assuming the same parametric input \cite{ECSW}. Lastly, it should be noted that the values reported in this paper are averaged over the total number of evaluations performed for each  case study.

\subsection{Toy Example: Two story building with nonlinear links} \label{toy}

This numerical case study simulates a 3D two story building modeled with nonlinear links. This is considered an educational example demonstrating the potential of the pROM approach before scaling it up on a 3D, real life structure. For the set up of this frame structure all beams and columns are of a rectangular cross section ($40cm \times 40cm$) and are considered made of reinforced concrete. The material parameters assumed are: Young modulus $E=21 GPa$, Poisson ratio $v=0.20$ and density $\rho=2400 \frac{kg}{m^3}$. The length of each frame is $l=7.5m$, the width is $w=5m$ and the height $h=3.2m$. Each node has six degrees of freedom, three translations and three rotations. 

\begin{figure*}[!hb]
        \centering
        \begin{subfigure}[b]{0.475\textwidth}   
            \centering 
            \includegraphics[width=1.0\textwidth]{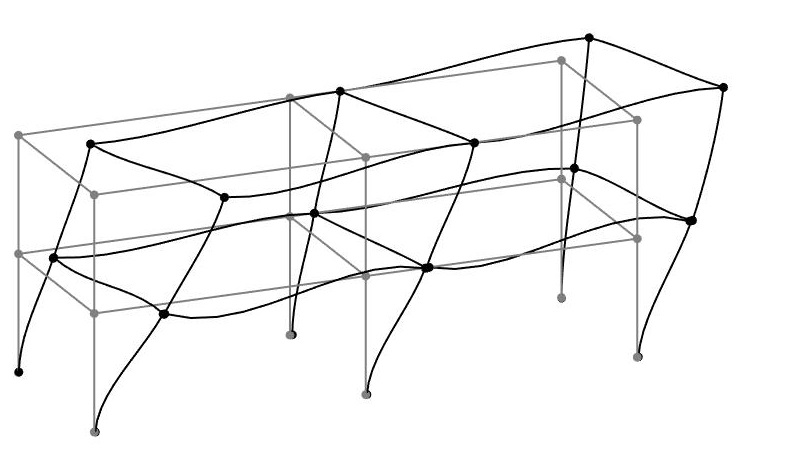}
            \caption{Geometrical Configuration of the problem.}
            \label{fig:boucwen}   
        \end{subfigure}
        \quad
        \begin{subfigure}[b]{0.475\textwidth}   
            \centering 
            \includegraphics[scale=1.0]{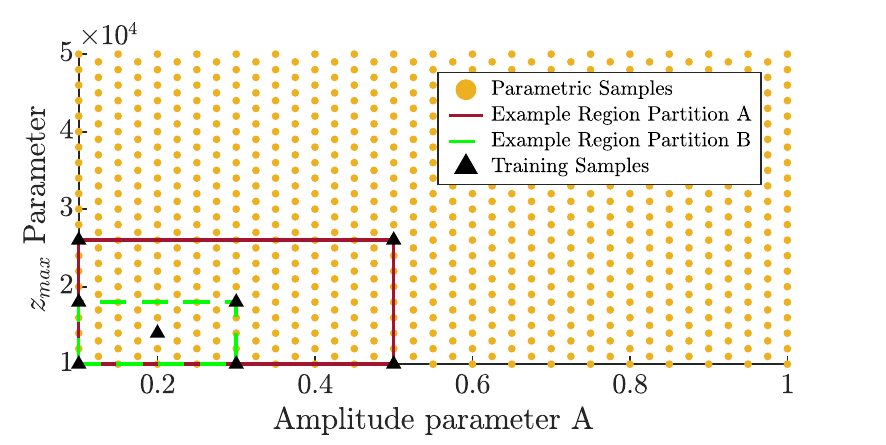}
            \caption{Domain sampling and partitioning approaches.} 
            \label{fig:region}  
        \end{subfigure}
        \caption{Two story building with nonlinear links. Geometrical configuration depicted in grey and example deformed state in black. Domain sampling and partitioning approaches examined are also visualized.} 
\end{figure*}

The structure is assumed fully constrained in the bottom nodes that represent the ground. Asymmetric nodal load excitation following a sinus form is modeled along with $4\%$ mass proportional damping. Each node of the structure is coupled with a hysteretic link modeled based on a nonlinear restoring force-displacement law. The same applies for the bottom element representing the ground. This way, all columns and beams are connected to the respective central and corner nodes through nonlinear links. This ensures that nonlinear links differentiate the global response sufficiently between parametric samples. The respective sketch of the geometrical configuration is presented with grey in Figure \ref{fig:boucwen}, including an example deformed geometry based on the modeled excitation presented in black. 

The Bouc Wen model is used for modeling the nonlinear restoring force contributions \cite{bw}. The respective Equation is:   

\begin{align}
\dot{\mathbf{z}} = \mathbf{A}\dot{x} - \beta |\dot{x}|\mathbf{z}|\mathbf{z}|^{n-1}-\gamma \dot{x}|\mathbf{z}|^{n} \nonumber
\\ 
\mathbf{z}_{max}= \Big( \frac{\mathbf{A}}{\beta+\gamma} \Big)^{\frac{1}{w}},
\label{boucwen}
\end{align}
where $x$ denotes the displacement of the link, $z$ the hysteretic parameter and $\mathbf{A}$ controls the hysteresis amplitude. Parameter $w$ along with the hysteretic parameters $\beta$ and $\gamma$ determine the basic shape of the hysteresis loop. Their absolute value is not of interest, but rather their sum/difference that may define a hardening or softening relationship. For this reason, in the context of this study the parametric dependency is modeled on the amplitude parameter $A$ and on $z_{max}$ excluding the already parametrized influence of $A$. The respective range of the parameters is $[0.10 - 1.00]$ for $A$ and $[1e04 - 5e04]$ for $z_{max}$. The sampling grid is depicted in Figure \ref{fig:region} as well. Based on the pROM approach of this paper described in section \ref{3}, the domain is partitioned manually in rectangular subdomains.

Here, two example partitions of the domain are discussed. The first case study is hereby referred to as Partition A. The extent of a subdomain is defined by a 0.4 unit variation on amplitude factor $A$ and a $1.6e04$ variation on the hysteretic parameter $z_{max}$. An example is demonstrated in Figure \ref{fig:region} with a purple line. This partition divides the domain in approximately 4 plus 4 overlapping subdomains. For example the presented subdomain in Figure \ref{fig:region} spans [$0.1$,$1.0e04$]-[$0.5$,$2.6e04$], whereas another subdomain would span [$0.5$,$1.0e04$]-[$0.1$,$2.6e04$], and so on. A finer sampling case study is also discussed, referred to as Partition B. This divides the domain in 20 subdomains plus 5 overlapping ones to account for the remaining range. An example is also presented in Figure \ref{fig:region} with a dashed green line and spans [$0.1$,$1.0e04$]-[$0.3$,$1.8e04$].

\begin{table}
\centering
	\caption[]
	{The performance of the pROM for the domain partition techniques presented in Figure \ref{fig:region}. The $\mathbb{R}\mathbb{E}_{\mathbf{rf}}$ error of Equation \eqref{errornorm} with respect to the restoring forces $rf$ is evaluated. The average and max error across the domain is presented along with the speed up factor. The pROM variants compared are described in Table \ref{tableref}.}
	\label{BoucWenTable}
	\begin{tabular} { P{4cm} P{2cm} P{2cm} P{2cm} P{2cm} }
	\toprule
		\multicolumn{5}{c}{Two story building with nonlinear links} \\ 
		\multicolumn{5}{c}{Partitioned Domain Span: [0.10,$1.0e04$]-[1.0,$5.0e04$]}\\
		\\
		 & Mean $\mathbb{R}\mathbb{E}_{\mathbf{u}}$ & Max $\mathbb{R}\mathbb{E}_{\mathbf{u}}$ & Mean $\mathbb{R}\mathbb{E}_{\mathbf{rf}}$ & Max $\mathbb{R}\mathbb{E}_{\mathbf{rf}}$ \\
		\hline
		\multicolumn{5}{c}{Partition A - Speed-up factor: 1.28}\\
		Global Basis & 3.18\% & 5.22\%&  4.84\% & 7.51\%\\
		Local Basis &  0.43\%& 0.68\% & 3.06\%  & 3.90\% \\
		Entries Interp. & 0.40\% & 0.67\%&  3.25\% & 3.73\%  \\
		Coefficients Interp. & 0.39\% & 0.67\%&  3.28\% & 3.69\% \\
		\multicolumn{5}{c}{Partition B - Speed-up factor: 1.31}\\
		Global Basis & 1.96\% & 5.91\%& 3.94\% & 8.28\%  \\
		Local Basis & 0.15\% & 0.33\% & 1.83\%   & 2.86\% \\
		Entries Interp. & 0.12\% & 0.30\%&  1.29\% & 2.86\% \\
		Coefficients Interp. & 0.12\% & 0.30\%&  1.30\% & 3.06\% \\
	\bottomrule
	\end{tabular}
\end{table}

The performance of the pROM variants presented in Table \ref{tableref} is summarized in Table \ref{BoucWenTable}. The average and maximum error measure for validation samples spanning across the parametric domain are presented. The accuracy is evaluated with respect to the response time history and the time history of the nonlinear terms, namely the restoring forces. By comparing the accuracy of the four pROM variants of Table \ref{tableref}, it seems that the use of a single projection basis on the Global Basis pROM of Table \ref{BoucWenTable} is not capable of approximating the underlying phenomena accurately enough. In addition, the proposed Coefficients Interpolation approach seems able to achieve a similar accuracy to the established element-wise Entries Interpolation and thus reproduce the high fidelity dynamic behavior and response.

Regarding computational reduction no hyper reduction technique is utilized for this case study. For this reason, the efficiency achieved is negligible due to the full evaluation of the nonlinear terms. The focus of this example lies on providing a 'proof of concept' for the approximation pROM strategy employed and for this reason the hyper reduction is only assembled on the 3D, actual scale example that follows.

\begin{figure*}[!hb]
        \centering
        \begin{subfigure}[b]{0.475\textwidth}
            \centering
            \includegraphics[scale=0.70]{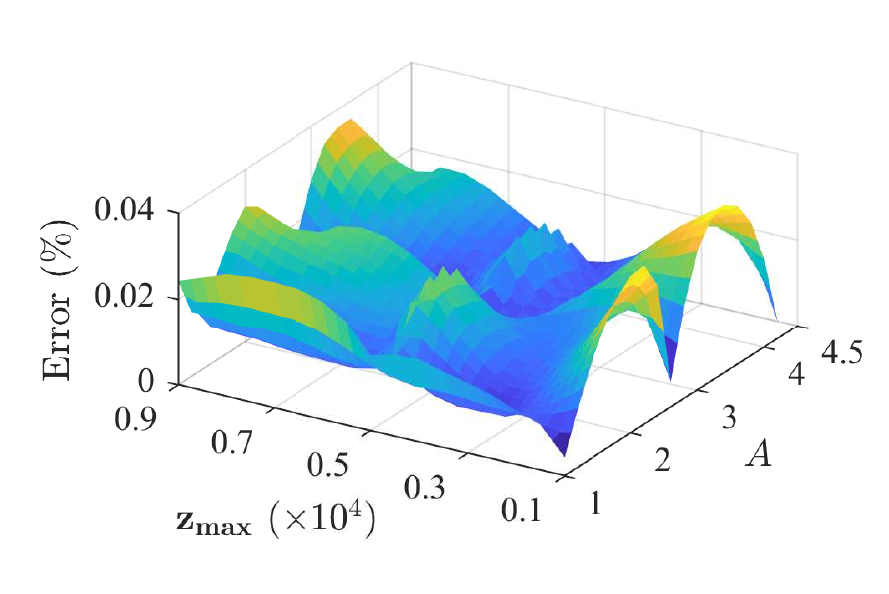}
            \caption[]%
            {Domain error plot for Partition A}    
            \label{fig:errorplot}
        \end{subfigure}
        \hfill
        \begin{subfigure}[b]{0.475\textwidth}  
            \hspace*{0.90cm} 
            \includegraphics[scale=0.60]{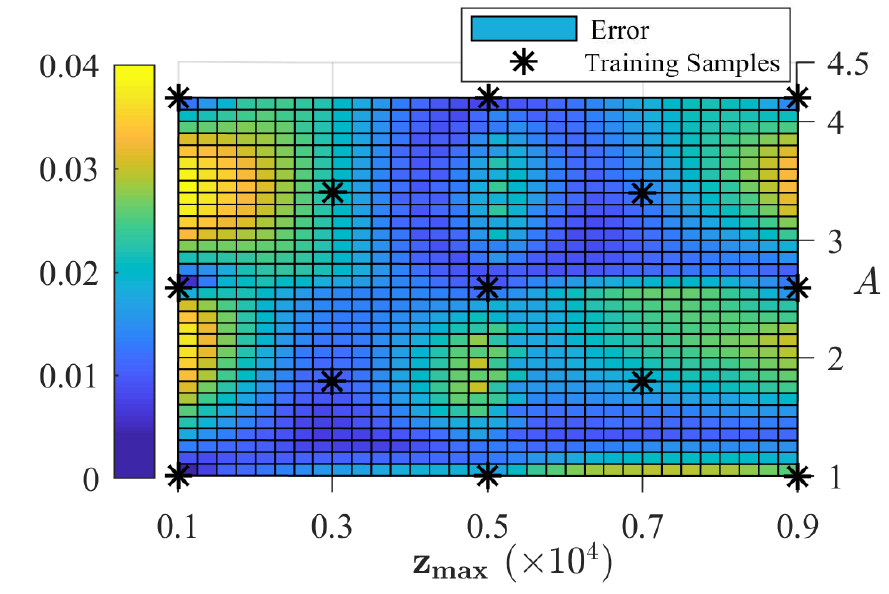}
            \caption[]%
            {Domain error in 2D for Partition A}    
            \label{fig:errorplotproject}
        \end{subfigure}
        \vskip\baselineskip
        \begin{subfigure}[b]{0.475\textwidth}   
            \centering 
            \includegraphics[scale=0.70]{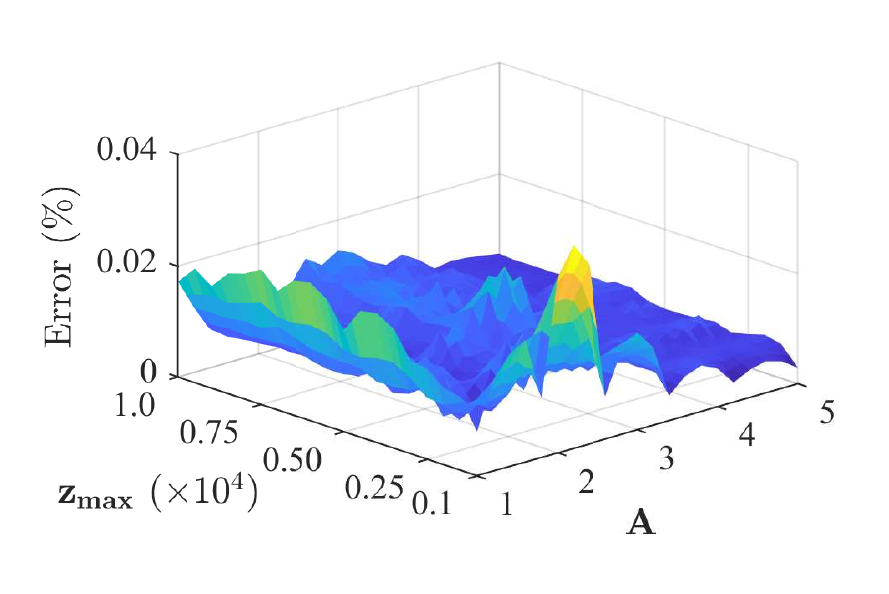}
            \caption[]%
            {Domain error plot for Partition B}    
            \label{fig:errorplotfine}
        \end{subfigure}
        \quad
        \begin{subfigure}[b]{0.475\textwidth}   
            \centering 
            \includegraphics[scale=0.60]{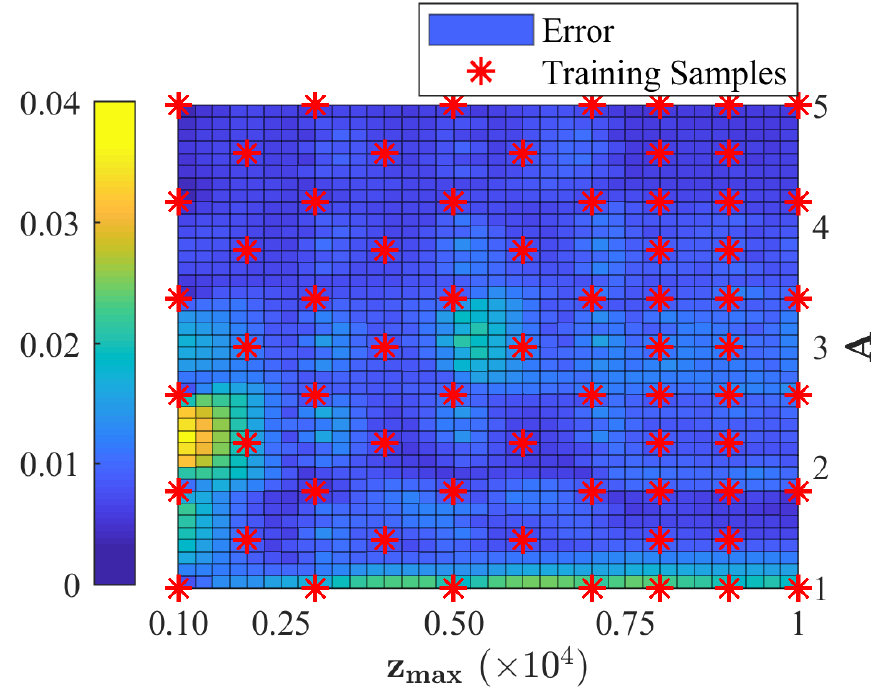}
            \caption[]%
            {Projection of domain error in 2D for Partition B}    
            \label{fig:errorplotprojectfine}
        \end{subfigure}
        \caption[]
        {Domain error plot for Bouc Wen model for Partition A and Partition B for the Coefficients Interpolation pROM (Table \ref{tableref}). The $\mathbb{R}\mathbb{E}_{\mathbf{rf}}$ error of Equation \eqref{errornorm} is evaluated with respect to the approximation of the restoring forces $rf$. A 3D and a 2D projection error plot are provided. } 
        \label{fig:ploterror}
\end{figure*}

The accuracy measure on displacements for the Coefficients Interpolation pROM (Table \ref{tableref}) proposed in this paper, is depicted in detail in Figure \ref{fig:ploterror}. A 3D and a 2D projection error plot are provided for each partition case. In Figure \ref{fig:errorplot} that corresponds to a coarse sampling of the domain, the pROM delivers an approximation of restoring forces with a relative error lower than $4\%$ in any case. As expected, Figures \ref{fig:errorplot} and \ref{fig:errorplotproject} indicate a lower error close to the training samples. The error also experiences a smooth arc transition between training samples, indicating an increase for increasing distance between validation and training samples. Figure \ref{fig:errorplotfine} presents the respective error for Partition B. Here the sampling is finer and therefore the relative error is substantially reduced across the whole domain. This implies an overall better approximation. In Figure \ref{fig:errorplotprojectfine} it is further observed that the error is minimized in the vicinity of the training samples, as expected.

\begin{figure*}
        \centering
        \begin{subfigure}[b]{0.475\textwidth}   
            \centering 
            \includegraphics[scale=0.90]{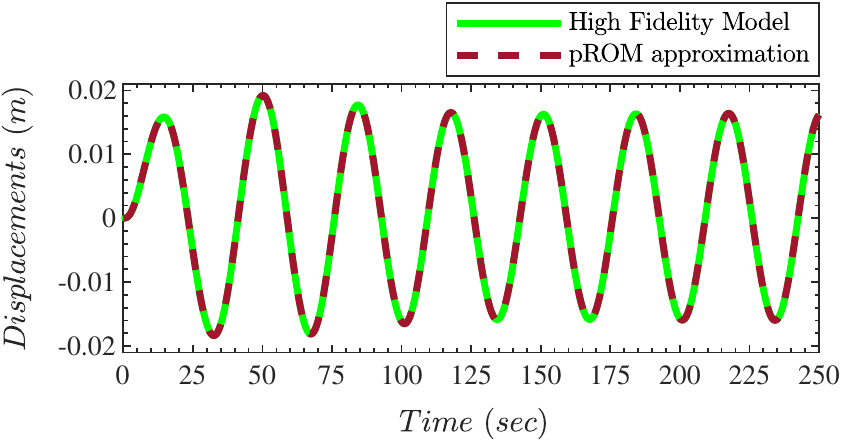}
            \caption{Response Time History of the model}
            \label{fig:BoucWenTH}   
        \end{subfigure}
        \quad
        \begin{subfigure}[b]{0.475\textwidth}   
            \centering
            \includegraphics[scale=0.90]{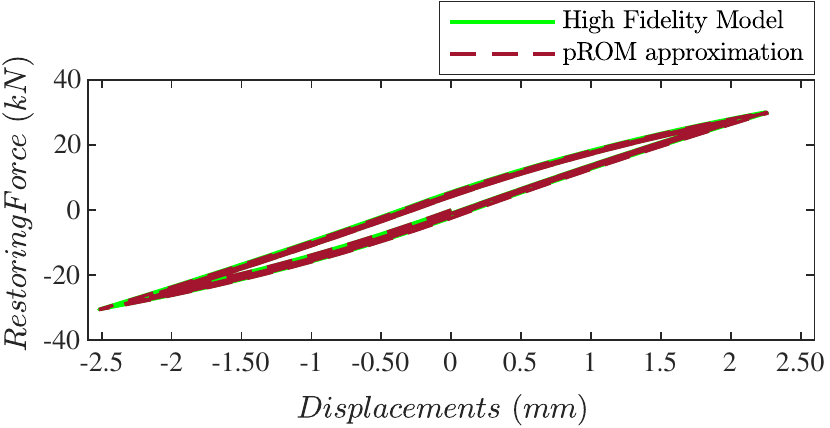}
            \caption{pROM approximation of steep hysteresis curve}
            \label{fig:bwprom1}   
        \end{subfigure}
        \vskip\baselineskip
        \begin{subfigure}[b]{0.475\textwidth}   
            \centering 
            \includegraphics[scale=0.90]{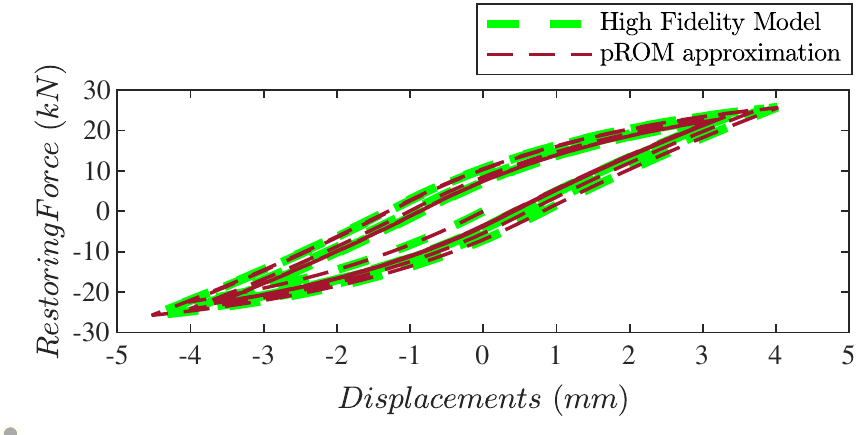}
            \caption{pROM approximation medium of hysteresis curve}
            \label{fig:bwprom2}   
        \end{subfigure}
         \quad
        \begin{subfigure}[b]{0.475\textwidth}   
            \centering 
            \includegraphics[scale=0.90]{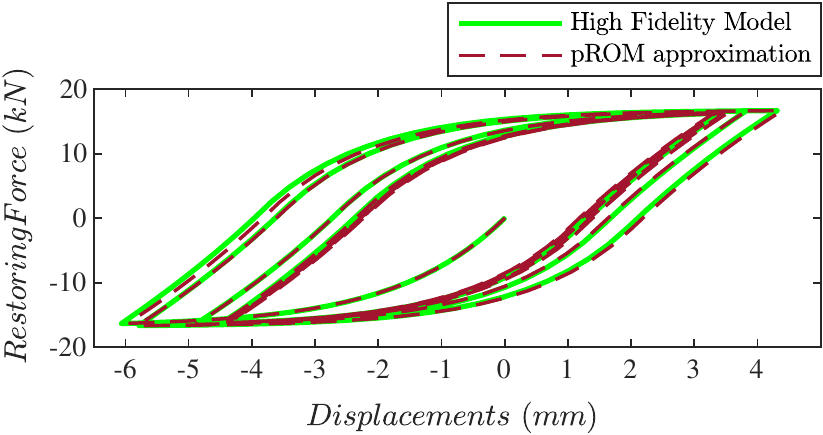}
            \caption[]%
            {pROM approximation of shallow hysteresis curve}    
            \label{fig:bwprom3}
        \end{subfigure}
        \caption{Accuracy of the proposed Coefficients Interpolation pROM (Table \ref{tableref}) for Partition B on capturing the variation of the hysteretic component. The approximation accuracy on a response time history and on three different shape of the hysteresis curve are depicted to demonstrate the potential of the method.} 
        \label{fig:bwresults}
\end{figure*}

In Figure \ref{fig:bwresults} the Coefficients Interpolation pROM (Table \ref{tableref}) proposed in this paper is evaluated for Partition B. The respective performance is validated in several samples to demonstrate the overall accuracy on estimating different shape and magnitude cases for the hysteresis curve. First, as illustrated in Figure \ref{fig:BoucWenTH}, the pROM seems able to reproduce the underlying HFM response as two response curves are practically indistinguishable. In Figures \ref{fig:bwprom1},\ref{fig:bwprom2},\ref{fig:bwprom3} the pROM accuracy in approximating the hysteresis curve is depicted. A steep case study is captured accurately in \ref{fig:bwprom1}, whereas a shallow case is approximated in Figure \ref{fig:bwprom3}. Figure \ref{fig:bwprom2} demonstrates a hysteresis curve lying somewhere in between a shallow and a steep shape. In all cases the respective error measure for the restoring forces in Equation \eqref{errornorm} lies below $3\%$. This indicates that the pROM is capable of accurately reproducing the underlying hysteresis phenomena and thus the respective response and dynamics for a range of parametric samples. 

Thus, the proposed pROM seems able to address case studies of nonlinear parametric dependency pertaining on properties of the system across a wide range of input. Efficiency gains employing hyper reduction are addressed in detail on the next case study. To this end, the framework presented here is scaled up to address a real life structure under earthquake excitation. For this example the dependency is formulated with respect to the amplitude and the frequency content of the earthquake spectrum, thus arguing about the potential of the proposed pROM on handling such case studies. 

\subsection{Wind Turbine Tower under parametric Earthquake Excitation} \label{tower}
\subsubsection{Numerical set up} \label{41}

The numerical case study is based on the simulated dynamic response of the NREL 5-MW baseline wind turbine tower \cite{Jonkman2009}, which is derived from a three-dimensional finite element model using shell elements. The circular cross-section of the tower is linearly tapered from the base, with diameter and thickness equal to 6m and 0.027m respectively, to the top, with 3.87m diameter and 0.019m thickness. The tower is made of steel with modulus equal to 210 GPa and material density of 7850 kg/m$^3$, which is slightly increased in order to account for additional structural components such as bolts, welds and flanges and further achieve a sufficient agreement of calculated vibration modes with the ones reported in \cite{Jonkman2009}. The tower is considered to be fully restrained at the base and the nonlinearity of the  model lies in the material constitutive law, which is characterized by isotropic von Mises plasticity.

The FE model is first assembled in a verified and established reference FEM software (ABAQUS \citep{abq}) using shell elements, and subsequently in MATLAB. Shell elements are integrated on the MATLAB framework based on the suggestions in \citet{Bathe1985}, while nonlinear behavior is modeled according to the assumptions and formulations in \citet{Bathe,Plastic}. The geometrical and material configuration of the tower, as well as the mesh properties are given in Table \ref{tablemat}. A constant thickness assumption throughout the tower is made for simplification purposes.

\begin{figure*}[!ht]
        \centering
        \begin{subfigure}[b]{0.475\textwidth}   
            \centering 
            \includegraphics[scale=0.35]{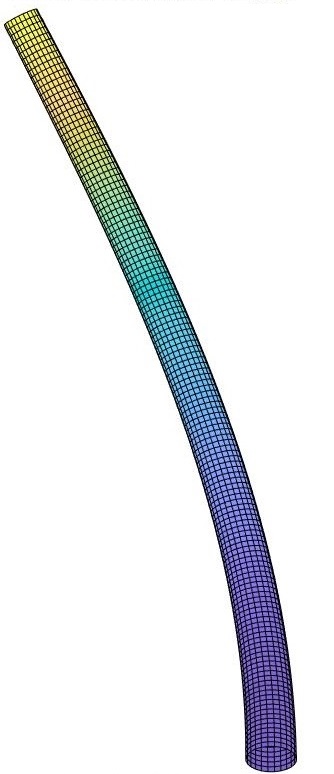}
            \caption{MATLAB FEM deformed based on 1st eigenmode}
            \label{fig:abqmodel}   
        \end{subfigure}
        \quad
        \begin{subfigure}[b]{0.475\textwidth}   
            \centering 
            \includegraphics[scale=0.35]{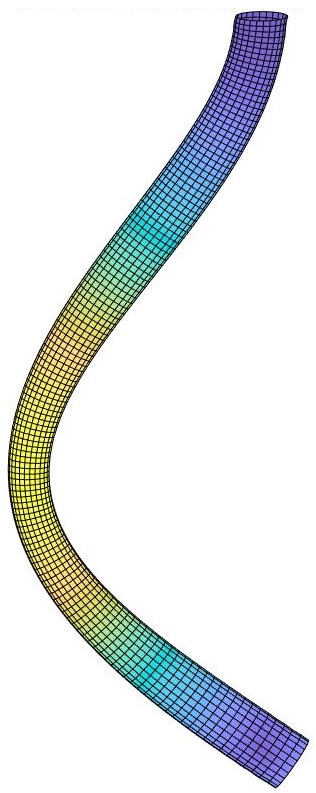}
            \caption{{MATLAB FEM deformed based on 5th eigenmode}}
            \label{fig:matmodel}   
        \end{subfigure}
        \caption{Finite Element Representation of the HFM}
        \label{fig:fem}
\end{figure*}

The ABAQUS FEM model is used as an independent reference representation for validation purposes. The full-order model representation assembled in MATLAB is used for the High Fidelity simulations and the respective parametric reduction framework. The respective HFM is presented in Figure \ref{fig:fem} in a deformed state to visualize the first and fifth eigenmode respectively. These eigenmodes are presented as the ones primarily excited due to the form of the earthquake orientation, namely only translation and bending modes are excited and not torsional or localized ones. 

\begin{table}[!hb]
\centering
	\caption{Mesh, Material and Geometrical configuration of the High Fidelity Model.}
	\label{tablemat}
	\begin{tabular} { p{4cm} p{3cm}  p{3cm} p{2cm} }
	\toprule
		\multicolumn{4}{c}{Wind turbine tower - High Fidelity Model}\\
		\hline
		Geometrical Properties & & Material Properties &\\
		Base Cross-Section & Radius: 3.00 m & Elastic Modulus, \texttt{E} & 210$e^{3}$ MPa\\
		Tip Cross-Section & Radius: 1.94 cm & 	Poisson Ration, $\nu$ & 0.30\\
		Height - Length & 87.61 m & Density,$\rho$  & 7850 $\frac{kg}{m^3}$\\
		Thickness & 19.00-27.00 mm & Yield Stress, \texttt{$f_{\mathrm{y}}$}& 435 MPa\\
		Mass of Wind turbine & 357 tons & Yield Surface & Von Mises\\
		\hline
		Mesh Properties & & &\\
		Mesh & \multicolumn{3}{l}{3264 4-node Linear Shell Elements}\\
		\hline
		Modal Properties & & &\\
		Eigenfrequencies & \multicolumn{3}{c}{\nth{1}=0.29 Hz, \nth{3}=2.97 Hz, \nth{5}=3.01 Hz, \nth{7}=4.52 Hz, \nth{9}=6.77 Hz}\\
		Damping & \multicolumn{3}{l}{Rayleigh Formulation}\\
		 & \multicolumn{3}{l}{2\% on \nth{1} \& \nth{5} eigenmodes}\\
	\bottomrule
	\end{tabular}
\end{table}

Regarding parametric dependency, the wind tower is subjected to various forcing inputs formulated as earthquake excitation. The forcing inputs are applied in the form of nodal loads, namely the acceleration of the amplitude spectrum is scaled with the respective mass and applied on all nodes of the HFM. The inputs are characterized by different temporal and spectral characteristics. This introduces the input uncertainty that is needed to explore the efficiency range of the reduction framework and assemble a reliable pROM across the whole domain of interest. Thus, the dependency of the ROM is represented in the excitation input of the model and the uncertainty parameters are the characteristics of the earthquake acceleration spectrum.

\begin{figure*}[!hb]
        \centering
        \begin{subfigure}[b]{0.485\textwidth}   
            \centering 
            \includegraphics[scale=0.90]{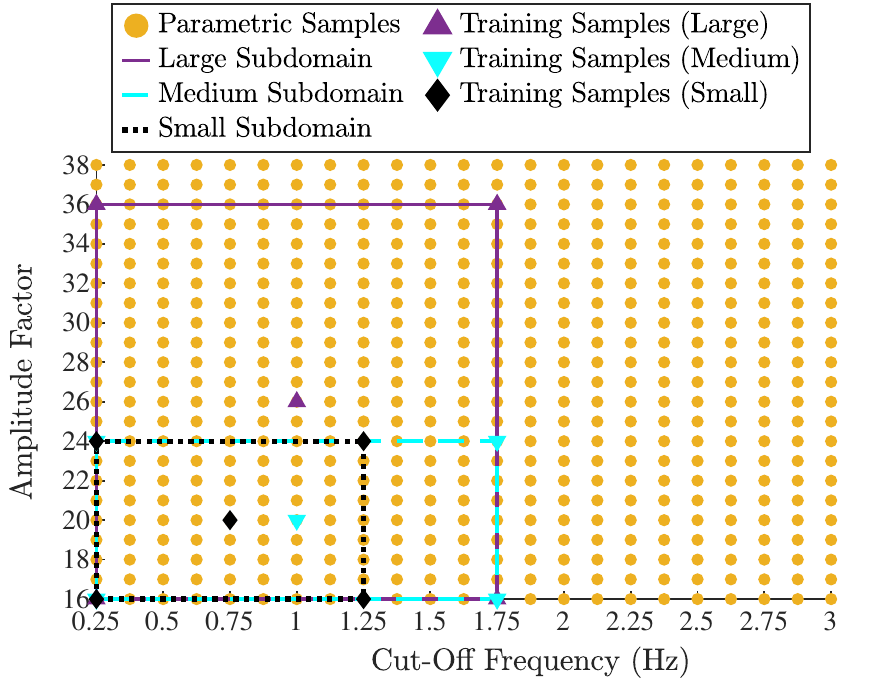}
            \caption{Parametric domain and partitioning regions.}
            \label{fig:grid} 
        \end{subfigure}
        \quad
        \begin{subfigure}[b]{0.485\textwidth}   
            \centering 
            \includegraphics[scale=0.90]{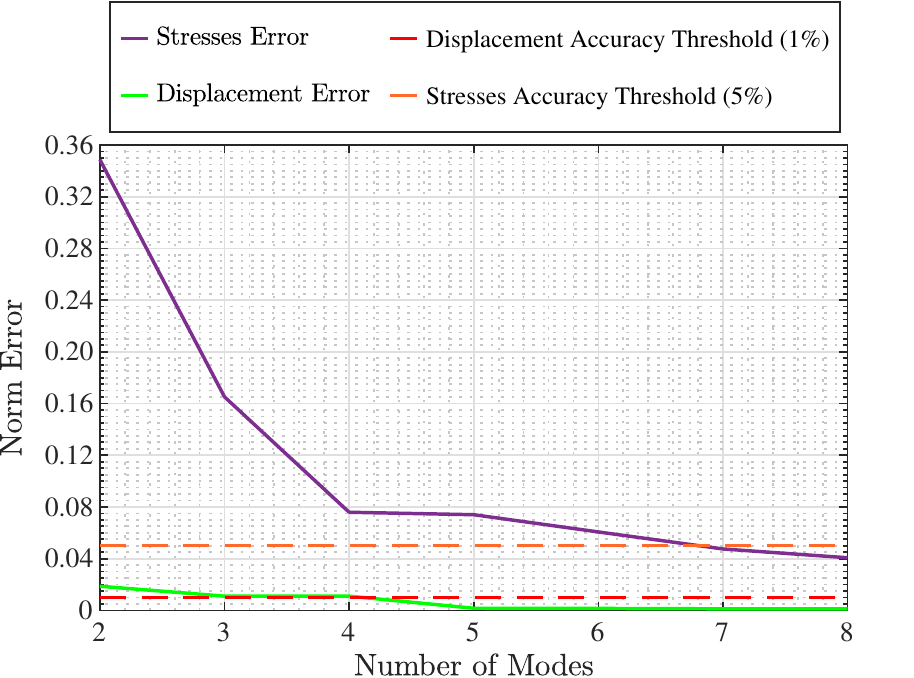}
            \caption{Selection of reduction order for the pROM.}
            \label{fig:basis} 
        \end{subfigure}
        \caption{Numerical set up considerations for the wind turbine tower case study. Reduction order selection based on accuracy thresholds and parametric domain of interest with respective partitions. The $\mathbb{R}\mathbb{E}_{\mathbf{u}}$ and $\mathbb{R}\mathbb{E}_{\mathbf{\sigma}}$ error of Equation \eqref{errornorm} for an example validation sample are presented. }
        \label{fig:setup}
\end{figure*}

In the context of this study, the earthquake excitation is approximated as a low-pass filtered random signal. Therefore two input parameters are required for parameterizing the input, namely the amplitude of the excitation and the cut-off frequency of the low-pass filter. 
For each sample, based on the respective parametric configuration, a white noise signal is passed through a low pass filter. The output of the filter is then scaled based on the amplitude parameter. The final processed signal represents the acceleration in an earthquake scenario and the HFM is simulated using this excitation. The duration of the excitation is 10 seconds, whereas the model is simulated for 30 seconds with a timestep of $\Delta t=0.01$. The simulations are carried out so as to succeed yielding of the tower, i.e., until a considerable amount of elements have entered the plastic region, corresponding to a distributed yielding domain.

The parametric dependency is herein expressed with respect to the temporal and spectral characteristics of the earthquake excitation. The respective amplitude parameter ranges from a factor of 16 to a factor of 38, whereas the respective frequency parameter from 0.25 up to 3Hz. The assumed frequency content excites at least the first three eigenmodes based on the set up parameters on Table \ref{tablemat}, an assumption sufficiently accurate to the real world response of the structure. Based on this range of parameters the respective parametric grid of interest is assembled. In Figure \ref{fig:grid} the domain of the parametric samples is depicted along with the different partition regions tested for the pROM training configuration. 

Moreover, the reduction order for the pROM was defined based on Figure \ref{fig:basis}. The $\mathbb{R}\mathbb{E}_{\mathbf{u}}$ and $\mathbb{R}\mathbb{E}_{\mathbf{\sigma}}$ error measures of Equation \eqref{errornorm} are presented with respect to the reduction order, namely the number of modes included in the projection basis. The thresholds for both error measures are defined based on the desired level of accuracy on the approximation. Based on Figure \ref{fig:basis}, a reduction order of 8 is chosen to ensure better approximation quality in both stresses and displacements while preserving the symmetrical properties of the problem.

\subsubsection{Results and Discussion}\label{42}
This subsection is devoted on validating the projected pROM approach of section \ref{3} numerically. The accuracy of the pROM is validated via comparison against the HFM. In addition, the pROM variants of Table \ref{tableref} are compared to study the potential and limitations of the proposed approach. The accuracy of the ROM approximation is studied with respect to displacements and stresses. Stresses are evaluated as they are fundamental quantities of interest in condition assessment, residual life estimation or any other kind of structural health monitoring application. Moreover, stresses represent a second order term of approximation and the respective error measure can be the determining factor of distinguishing the superiority or effectiveness of the employed approximation technique.

Regarding the sampling scheme, rectangular subdomains of the grid are defined and interpolation is performed inside. The extent and the location of the subdomains vary to validate the consistency and limitations of the proposed approach. The rectangular shape of the subdomains is chosen for simplification as the focus of the paper lies elsewhere. The limitations introduced by such simplifying assumptions are addressed in section \ref{5}. 

The example partitions in subdomains used here to validate the accuracy of the pROM are depicted in Figure \ref{fig:grid}. Three subdomains are examined in this paper, namely a small, a medium and a large one, depending on their size and range. For each rectangular region the training samples are located on the four edges and the centroid that also serves as reference point for the projection to the tangent plane. These are depicted with triangles and diamonds in Figure \ref{fig:grid}. Training and validation samples are hereby referred based on the parametric configuration, namely $[$Cut-Off Frequency, Amplitude Factor$]$. For instance, in Figure \ref{fig:grid} the large subdomain example presented spans [0.25,16]-[1.75,36] whereas another example in that case could span [1.75,18]-[3.00,38] and so on.

As an initial approach, large subdomains are formed. Local bases interpolation relies on the fact that the dynamics on the vicinity of each parameter sample can be spanned by a local subspace. Therefore, since the principal independent components of the response time histories of neighboring parameter samples are assumed to span the same subspace, the respective training samples should be located 'close' enough. By initially selecting large subdomains we aim to validate this argument by negation. The respective large subdomain example is depicted in Figure \ref{fig:grid} with a purple continuous line. 

Table \ref{OneRectangle} summarizes the accuracy performance of the pROM for large subdomains. The  $\mathbb{R}\mathbb{E}_{\mathbf{\sigma}}$ error measure of Equation \eqref{errornorm} is presented in detail for example validation snapshots of the example subdomain of Figure \ref{fig:grid} spanning [0.25,16]-[1.75,36] along with an average and maximum error value for equivalent subdomains spanning the rest of the domain. 

\begin{table}
\centering
	\caption{pROM results for large subdomains. The subdomain definition is based on the respective notation on Figure \ref{fig:grid}. The $\mathbb{R}\mathbb{E}_{\mathbf{\sigma}}$ error of Equation \eqref{errornorm} is evaluated. Detailed errors for validation samples in an example subdomain are presented, along with the average and max error of equivalent subdomains spanning the rest of the domain. The approaches compared are based on the variant pROMs of Table \ref{tableref}.}
	\label{OneRectangle}
		\begin{tabular} { P{4cm} P{1.50cm} P{1.50cm} P{1.50cm} P{1.50cm} P{1.75cm} P{1.50cm} }
		\toprule
		\multicolumn{7}{c}{Large Subdomains Scenario $\mathbb{R}\mathbb{E}_{\mathbf{\sigma}}$ error metric}\\ 
		\multicolumn{7}{c}{Span of Example Subdomain: [0.25,16]-[1.75,36],
		Partitioned Domain Span: [0.25,16]-[3.00,38]}\\
		\\
		 & Validation Sample [0.62,31] & Validation Sample [0.62,21] & Validation Sample [1.37,21] & Validation Sample [1.37,31]& Average $\mathbb{R}\mathbb{E}_{\mathbf{\sigma}}$ error (Domain) & Max $\mathbb{R}\mathbb{E}_{\mathbf{\sigma}}$ error (Domain)  \\
		\hline
		Global Basis& 12.99\% & 11.12\%& 7.08\% & 7.67\% & 9.72\% & 12.99\%\\
		Local Basis & 5.65\%& 7.16\%& 4.91\% & 4.75\% & 5.65\%& 7.16\%\\
		Entries Interp. & 10.36\% & 11.42\% & 10.03\% & 7.99\% & 9.92\%& 11.42\%\\
		Coefficients Interp.& 10.47\%& 11.35\%& 10.01\%&8.11\% & 9.96\% & 11.35\%\\
	\bottomrule
	\end{tabular}
\end{table}

\begin{figure*}[!ht]
        \begin{subfigure}[b]{0.475\textwidth}
            \includegraphics[scale=0.80]{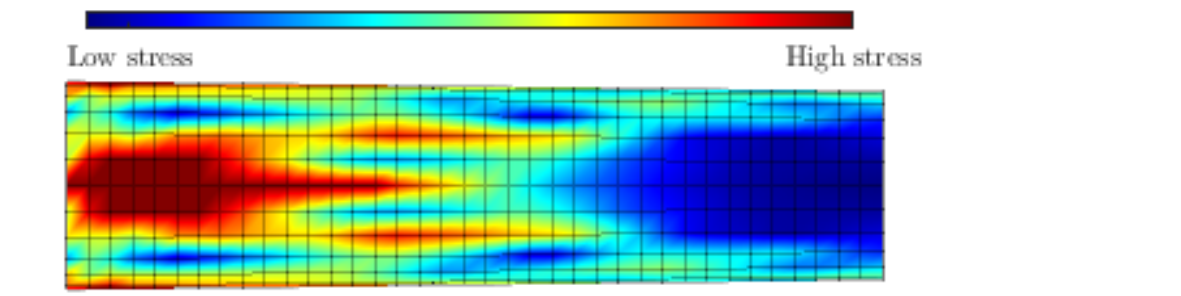}
            \caption[]%
            {Stress visualization for the High Fidelity Model}    
            \label{fig:wide}
        \end{subfigure}
        \hfill
        \begin{subfigure}[b]{0.475\textwidth}  
            \includegraphics[scale=0.85]{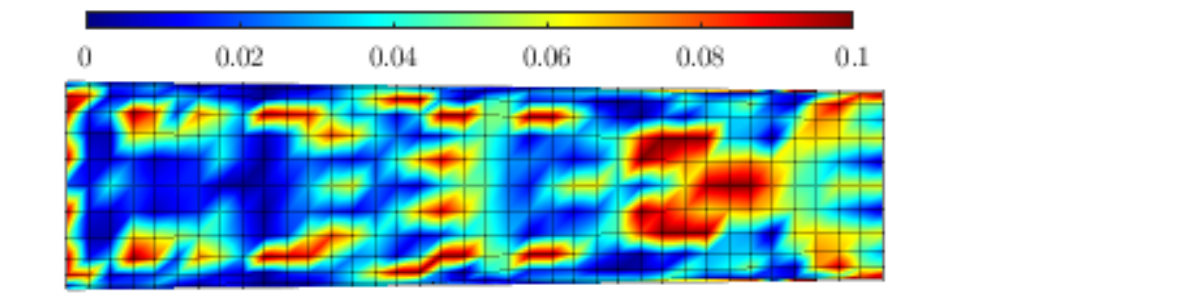}
            \caption[]%
            {Local Basis approximation error}    
            \label{fig:wideG}
        \end{subfigure}
        \vskip\baselineskip
        \begin{subfigure}[b]{0.475\textwidth}   
            \hspace*{-0.14cm}
            \includegraphics[scale=0.85]{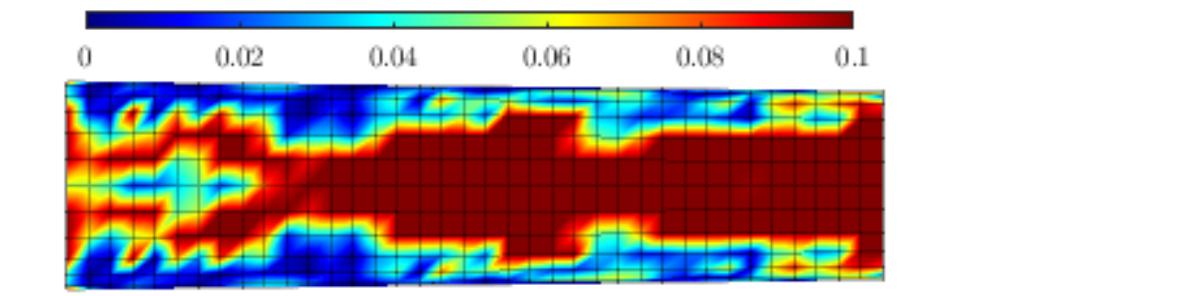}
            \caption[]%
            {Coefficients Interpolation approximation error}    
            \label{fig:wideR}
        \end{subfigure}
        \quad
        \begin{subfigure}[b]{0.475\textwidth}   
            \hspace*{0.28cm}
            \includegraphics[scale=0.85]{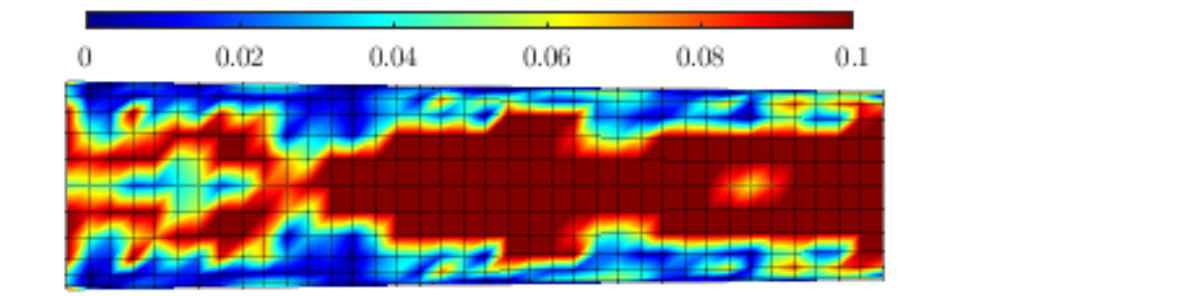}
            \caption[]%
            {Entries Interpolation approximation error}    
            \label{fig:wideEl}
        \end{subfigure}
        \caption{Stress approximation accuracy of the pROM for large subdomains. One third of the height of the wind tower is visualized. Stresses are visualized for the [0.62,31] validation sample of Table \ref{OneRectangle}. The effective von Mises stress \cite{Plastic} is depicted using node averaging and the $\mathbb{R}\mathbb{E}_{\mathbf{\sigma}}$ error is presented from Equation \eqref{errornorm}. The pROM are referred using the notation of Table \ref{tableref}.}  
        \label{fig:wideerror}
\end{figure*}

It can be inferred from Table \ref{OneRectangle} that the Local Basis approach delivers more accurate results in this case whereas the Global Basis performance implies that this variant might not be effective on nonlinear dynamic problems. Complementary to Table \ref{OneRectangle}, Figure \ref{fig:wideerror} depicts the stress contours of the wind turbine tower. A 2D projection for one third of the height of the wind turbine tower is presented, namely the yielding domain. An equivalent visualization can be observed for the rest of the circular cross section. Figure \ref{fig:wide} represents the results for the HFM, while Figure \ref{fig:wideG},\ref{fig:wideR} and \ref{fig:wideEl} the respective approximation $\mathbb{R}\mathbb{E}_{\mathbf{\sigma}}$ error of Equation \eqref{errornorm} using the Local Basis, the Coefficients Interpolation and the Entries Interpolation variants respectively. These findings suggest that the Local Basis pROM performs better in reproducing the stress contours of the High Fidelity Model for this sampling scheme.

For the poor performance in the case of interpolation based pROMs, the coarse sampling has to be considered. Specifically, after the tangent space projection mapping, the training samples might 'lie far' from the tangent plane drawn on the reference point on the centroid of the subdomain. For this reason, both interpolation schemes may deliver a subspace unable to approximate the desired solutions accurately. To this end, the domain is partitioned in finer subdomains. The respective example is depicted in Figure \ref{fig:grid} with a cyan dashed line. This case study is referred to as medium subdomains, with 'medium' referring to size comparison between the subdomains cases of Figure \ref{fig:grid}.

\begin{figure*}[!ht]
        \begin{subfigure}[b]{0.47\textwidth}
            \includegraphics[scale=0.80]{media/Stresses3162.pdf}
            \caption[]%
            {Stress visualization for the High Fidelity Model}    
            \label{fig:medium}
        \end{subfigure}
        \hfill
        \begin{subfigure}[b]{0.475\textwidth}  
            \includegraphics[scale=0.85]{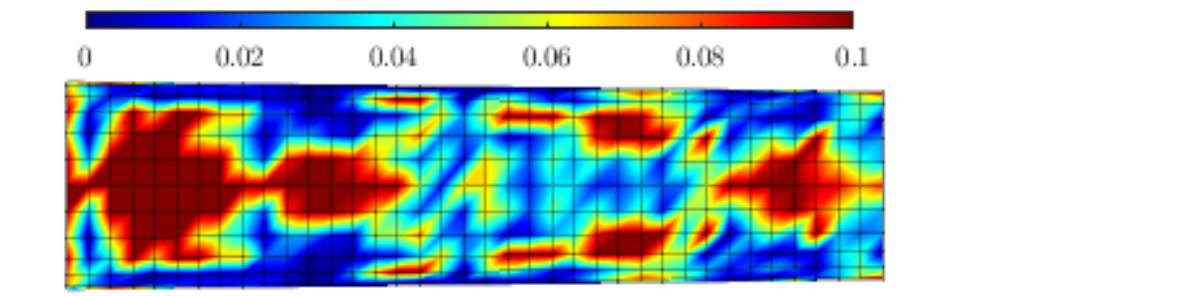}
            \caption[]%
            {Local Basis approximation error}    
            \label{fig:mediumG}
        \end{subfigure}
        \vskip\baselineskip
        \begin{subfigure}[b]{0.47\textwidth}   
            \hspace*{-0.14cm}
            \includegraphics[scale=0.85]{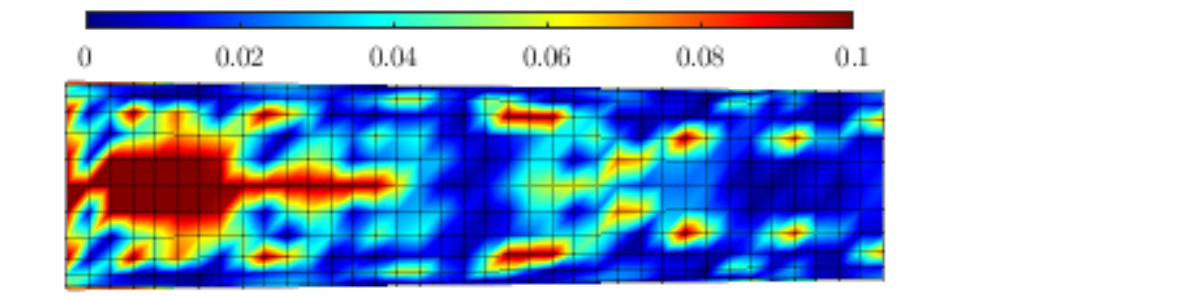}
            \caption[]%
            {Coefficients Interpolation approximation error}    
            \label{fig:mediumR}
        \end{subfigure}
        \quad
        \begin{subfigure}[b]{0.475\textwidth}   
            \hspace*{0.35cm}
            \includegraphics[scale=0.85]{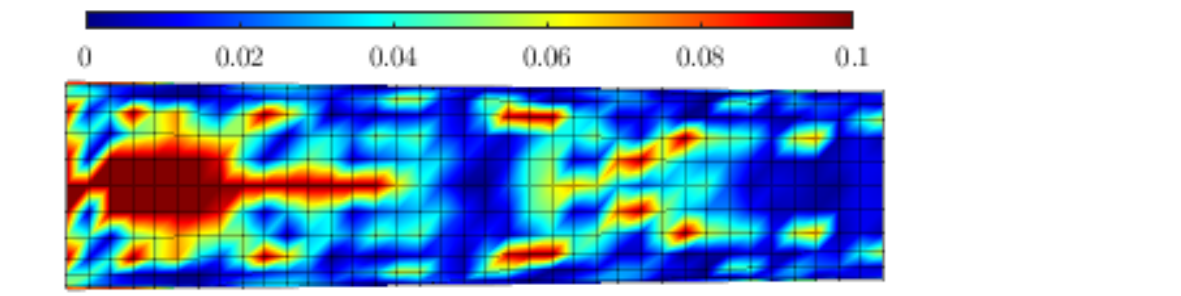}
            \caption[]%
            {Entries Interpolation approximation error}    
            \label{fig:mediumEl}
        \end{subfigure}
        \caption[ Stress approximation for a wide localized region.]
        {Stress approximation accuracy of the pROM for medium range regions. One third of the height of the wind tower is visualized. Stresses are visualized for the [0.62,31] validation sample of Table \ref{CompositeRectangle}.The effective von Mises stress \cite{Plastic} is depicted using node averaging and the $\mathbb{R}\mathbb{E}_{\mathbf{\sigma}}$ error is presented from Equation \eqref{errornorm}. The pROMs are referred using the notation of Table \ref{tableref}.}  
        \label{fig:mediumerror}
\end{figure*}

The domain of interest is divided in four subdomains. Based on Figure \ref{fig:grid}, these regions span: [$0.25$,$16$]-[$1.75$,$24$], [$0.25$,$28$]-[$1.75$,$36$], [$1.50$,$16$]-[$3.00$,$24$] and [$1.50$,$28$]-[$3.00$,$36$]. The performance of the pROM variants of Table \ref{tableref} is evaluated and the results are summarized on Table \ref{CompositeRectangle}. Validation is further provided on the samples defined in Table \ref{OneRectangle} for comparison purposes. In addition, Figure \ref{fig:mediumerror} depicts the stress contours of the wind turbine tower, similar to the large domain case study of Figure \ref{fig:wideerror}. The performance visualization on these figures is provided for the exact same validation sample for comparison purposes.

\begin{table}
	\caption{pROM results for medium subdomains. The subdomain definition is based on the respective notation on Figure \ref{fig:grid}. The $\mathbb{R}\mathbb{E}_{\mathbf{\sigma}}$ error measure of Equation \eqref{errornorm} is evaluated. Detailed errors for validation samples in an example subdomain are presented, along with the average and max error of equivalent regions spanning the rest of the domain. The approaches compared are based on the variant pROMs of Table \ref{tableref}.}
	\label{CompositeRectangle}
		\begin{tabular} { P{4cm} P{1.50cm} P{1.50cm} P{1.50cm} P{1.50cm} P{1.75cm} P{1.50cm} }
		\toprule
		\multicolumn{7}{c}{Medium Subdomains Scenario $\mathbb{R}\mathbb{E}_{\mathbf{\sigma}}$ error metric}\\ 
		\multicolumn{7}{c}{Span of Example Subdomain: [0.25,16]-[1.75,24] and [0.25,28]-[1.75,36]}\\
		\multicolumn{7}{c}{Partitioned Domain Span: [0.25,16]-[3.00,36]}\\
		\\
		  & Validation Sample [0.62,31] & Validation Sample [0.62,21] & Validation Sample [1.37,21] & Validation Sample [1.37,31]& Average $\mathbb{R}\mathbb{E}_{\mathbf{\sigma}}$ error (Domain) & Max $\mathbb{R}\mathbb{E}_{\mathbf{\sigma}}$ error (Domain)  \\
		\hline
		Global Basis& 12.99\% & 11.12\%& 7.08\% & 7.67\% &8.43\% & 13.44\%\\
		Local Basis & 6.80\%& 6.98\%& 6.64\% & 3.89\% & 5.64\%& 8.46\%\\
		Entries Interp. & 4.08\% & 6.01\% & 5.48\% & 2.96\% & 4.02\%& 6.01\% \\
		Coefficients Interp.& 4.07\%& 6.10\%& 5.51\%&2.78\% & 3.98\%& 6.10\% \\
	\bottomrule
	\end{tabular}
\end{table}

Contrary to the results of Table \ref{OneRectangle}, Table \ref{CompositeRectangle} outlines interpolation variants that are more accurate. The average error across the examined regions lies around $4\%$ with a maximum value around $6\%$. Both interpolation approaches deliver similar results and are at any case around $1\%$ more accurate than the Local Basis variant. This argument is validated by comparing the stress contours in Figures \ref{fig:wideerror} and \ref{fig:mediumerror} respectively. The $\mathbb{R}\mathbb{E}_{\mathbf{\sigma}}$ error of the interpolation approaches on medium subdomains depicted in Figures \ref{fig:mediumR},\ref{fig:mediumEl}, reduces compared to Figures \ref{fig:wideR}, \ref{fig:wideEl}. This suggests that the pROM delivers a more accurate approximation in comparison to large subdomains. Moreover, the error on the interpolation approaches in Figures \ref{fig:mediumR},\ref{fig:mediumEl} is lower compared to the Local Basis approximation in Figures \ref{fig:wideG} and \ref{fig:mediumG}. 

All in all, the interpolation variants of Table \ref{tableref} seem accurate in approximating the stress state of the wind turbine tower. This validates our hypothesis that the extents of the defined subdomains play an important role to the accuracy of the pROM strategy. This seems to be mainly attributed to the projected snapshots lying on or close to the tangent space of the reference point as described previously. However, the error measure on the interpolation approaches is still considered high. Therefore a subsequent sampling refinement is to be made, to further examine the potentials of the pROM variants under consideration.

\begin{table}[!ht]
	\caption{pROM results for small subdomains. The subdomain definition is based on the respective notation on Figure \ref{fig:grid}. The $\mathbb{R}\mathbb{E}_{\mathbf{\sigma}}$ error measure of Equation \eqref{errornorm} is evaluated. Detailed errors for validation samples in an example subdomain are presented, along with the average and max error of equivalent regions spanning the rest of the domain. The approaches compared are based on the variant pROMs of Table \ref{tableref}.}
	\label{NarrowRectangle}
		\begin{tabular} { P{4cm} P{1.50cm} P{1.50cm} P{1.50cm} P{1.50cm} P{1.75cm} P{1.50cm} }
		\toprule
		\multicolumn{7}{c}{Small Subdomains Scenario $\mathbb{R}\mathbb{E}_{\mathbf{\sigma}}$ error metric}\\ 
		\multicolumn{7}{c}{Span of Example Subdomain: [1.00,30]-[2.00,38]}\\
		\multicolumn{7}{c}{Partitioned Domain Span: [0.25,16]-[3.00,38]}\\
		\\
		& Validation Sample [1.25,32] & Validation Sample [1.75,32] & Validation Sample [1.25,36] & Validation Sample [1.75,36]& Average $\mathbb{R}\mathbb{E}_{\mathbf{\sigma}}$ error (Domain) & Max $\mathbb{R}\mathbb{E}_{\mathbf{\sigma}}$ error (Domain)  \\
		\hline
		Global Basis& 8.16\% & 6.68\%& 8.68\% & 7.53\% &9.92\% & 14.48\%\\
		Local Basis & 4.90\%& 4.68\%& 4.23\% & 3.40\% & 4.45\%& 7.24\%\\
		Entries Interp. & 2.41\% & 2.43\% & 2.90\% & 2.20\% & 3.57\%& 6.30\% \\
		Coefficients Interp.& 2.31\%& 2.37\%& 2.84\%&2.10\% & 3.59\%& 6.24\% \\
		\bottomrule
	\end{tabular}
\end{table}

\begin{figure*}[!ht]
        \centering
        \begin{subfigure}[b]{0.475\textwidth}
            \centering
            \includegraphics[scale=0.90]{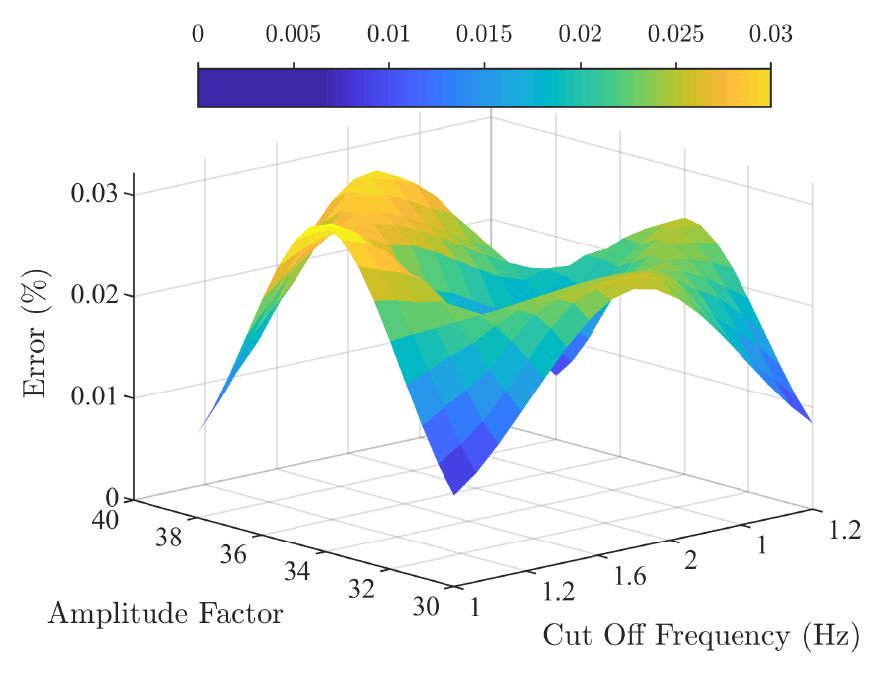}
            \caption[]%
            {Subdomain error plot of Table \ref{NarrowRectangle}}
            \label{fig:errorplotnarrow}
        \end{subfigure}
        \hfill
        \begin{subfigure}[b]{0.475\textwidth}  
            \centering 
            \includegraphics[scale=0.85]{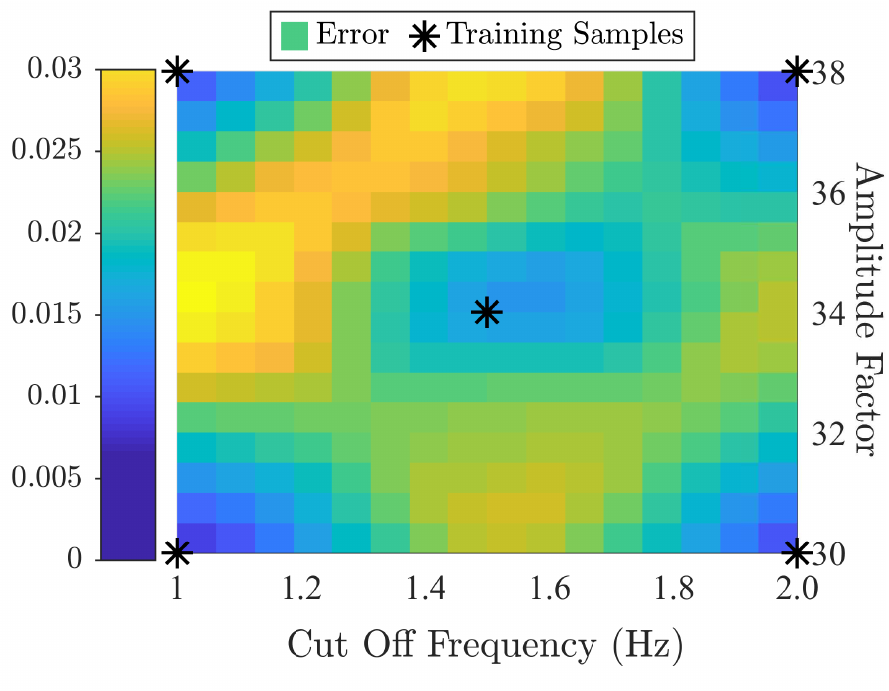}
            \caption[]%
            {Projection of subdomain error of Table \ref{NarrowRectangle}}    
            \label{fig:errornarrowproject}
        \end{subfigure}
        \caption{Subdomain error plot of the wind tower model for the small example subdomain of Table \ref{NarrowRectangle}. A 3D and a 2D projection error pot are provided. The $\mathbb{R}\mathbb{E}_{\mathbf{\sigma}}$ error measure of Equation \eqref{errornorm} is evaluated.}
        \label{fig:narrow}
\end{figure*}

In this case, the refined subdomains span between [$0.25$,$16$] and [1.25,24] and so on and are referred to as small subdomains. An example is demonstrated in Figure \ref{fig:grid} with a dotted black line. The respective results are summarized in Table \ref{NarrowRectangle} and Figure \ref{fig:narrow}. Since the domain partitioning changes, the example validation samples change as well. The validation samples presented in detail in Tables \ref{OneRectangle}, \ref{CompositeRectangle} are located on the middle of the diagonal distance between a training sample on the edge and the reference point on the centroid of the respective subdomain. The same principle is applied in Table \ref{NarrowRectangle} to demonstrate a fair comparison.  

Similar to Table \ref{CompositeRectangle}, the local bases interpolation pROMs seem more efficient in terms of accurately predicting the stress state of the model for the small subdomains of Table \ref{NarrowRectangle}. The Local Basis approach yields a better performance compared to Table \ref{CompositeRectangle}, however the respective error metrics remain higher. Comparing the interpolation approaches, it seems that the suggested approach of section \ref{3} performs marginally better than local bases direct interpolation. Taking into account that the Entries Interpolation technique also depends on the dimension $n$ of the HFM, this implies that the proposed approach of this study may produce a marginally better pROM in total. However, this marginal difference seems to be dependent on sampling and subdomain size parameters and this suggests that additional research and argumentation is needed for a definite conclusion with general application. All in all though, the suggested approach performs at least with the same accuracy with the respective local bases entry interpolation across the validation cases of this paper. 

The accuracy of the proposed pROM across the example subdomain of Table \ref{NarrowRectangle} is visualized in Figure \ref{fig:narrow}. A 3D and a 2D projection error plot with respect to stresses approximation are provided. Figure \ref{fig:errorplotnarrow} demonstrates that the pROM delivers a decent approximation regarding the $\mathbb{R}\mathbb{E}_{\mathbf{\sigma}}$ error with values lower than $3\%$ in any case. Figure \ref{fig:errornarrowproject} presents the same error plot in a 2D domain, indicating the positions of the training samples as well. 

\begin{figure*}[!hb]
        \begin{subfigure}[b]{0.475\textwidth}   
            \centering 
            \includegraphics[scale=0.85]{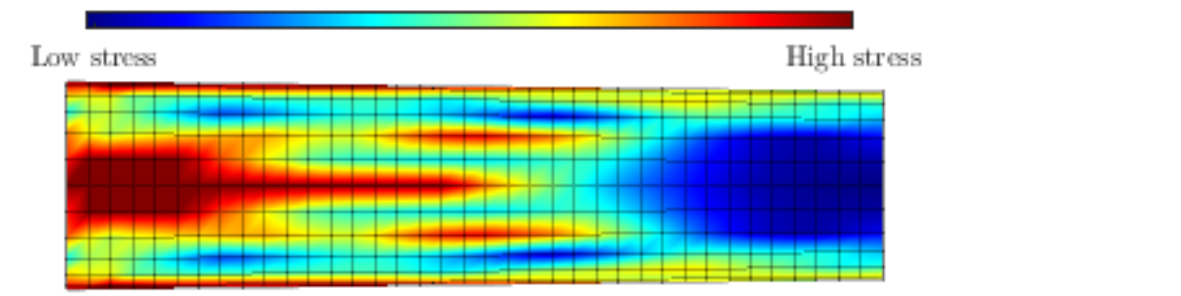}
            \caption[]%
            {Stress visualization for the High Fidelity Model}    
            \label{fig:stressHFM}
        \end{subfigure}
        \quad
        \begin{subfigure}[b]{0.475\textwidth}   
            \centering 
            \includegraphics[scale=0.85]{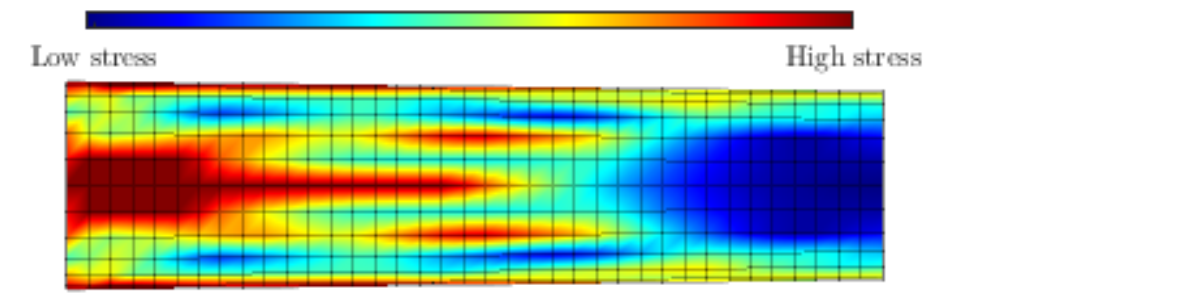}
            \caption[]%
            {Stress visualization for the pROM approximation} 
            \label{fig:stressROM}
        \end{subfigure}
        \vskip \baselineskip
        \begin{subfigure}{1.0\textwidth}   
            \hspace{3.5cm}
            \includegraphics[scale=1.00]{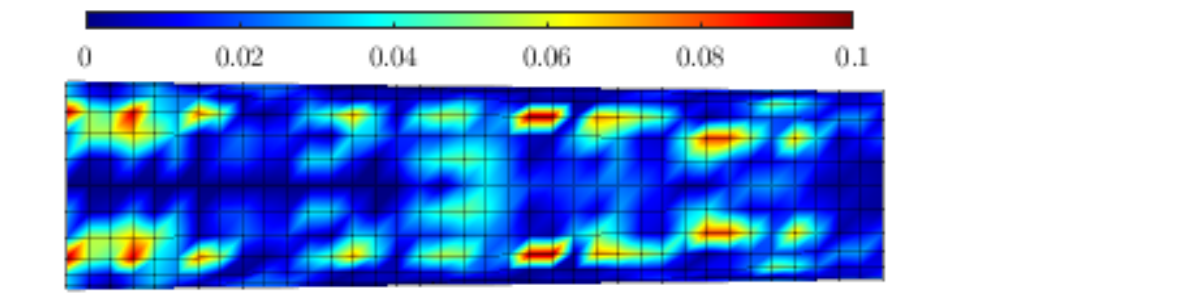}
            \caption[]%
            {Stress error visualization}
            \label{fig:stresserror}
        \end{subfigure}
        \caption{Stress visualization for the HFM and the pROM proposed.  One third of the height of the wind tower is visualized. Stresses are visualized for the [1.25,36] validation sample of Table \ref{NarrowRectangle}. The $\mathbb{R}\mathbb{E}_{\mathbf{\sigma}}$ error of Equation \eqref{errornorm} regarding the effective von Mises stress \cite{Plastic} is depicted using node averaging.} 
        \label{fig:stresses}
\end{figure*}

In addition, Figure \ref{fig:stresses} depicts the stress visualization and approximation of the yielding domain of the wind turbine tower in a qualitative manner. The results represent the validation sample [1.25,36] evaluated in Table \ref{NarrowRectangle} as well. Figure \ref{fig:stressHFM} presents the High Fidelity stress contours of the model and Figure \ref{fig:stressROM} the pROM approximation. Figure \ref{fig:stresserror} shows a visualization of the relative $\mathbb{R}\mathbb{E}_{\mathbf{\sigma}}$ error, offering a closer look on the quality of the approximation. Comparing these figures, it seems that the resulting pROM approximation delivers a sufficiently accurate estimation in reproducing the stress contours of the HFM. The qualitative visualization between the HFM and the pROM is almost identical, whereas the relative error in Figure \ref{fig:stresserror} assumes values lower than $4\%$ in the middle region, where the stress reaches its peak in Figures \ref{fig:stressHFM} and \ref{fig:stressROM}. 

In summary, the proposed pROM approach has been demonstrated to accurately reproduce the response of a nonlinear, large scale system with parametric dependencies pertaining to the characteristics of the earthquake excitation. However, the computational load required for these evaluations has not yet been addressed. The following subsection discuss this aspect. 

\subsubsection{ Hyper-Reduced pROM}

The hyper reduction strategy described in section \ref{sechyper} has been implemented here to assemble a hyper pROM, able to deliver substantial speed up on computational time. Equations \eqref{xi2},\eqref{HPmatrixes},\eqref{xi} are responsible for deriving the respective reduced mesh and weight coefficients for a sparse evaluation of the nonlinear terms of the problem.

For all evaluations an implicit Newmark integration scheme is implemented. As already described, the time discretization remains identical, namely a constant time-step is used, and a single core is utilized for each realization of either the HFM or the pROM. This establishes an objective comparison in terms of computation load. 

The ECSW hyper reduction method implemented here is based on producing a reduced mesh for the evaluation of the nonlinear terms \cite{ECSW}. For this reason, the achieved sparsity on the number of yielding elements being evaluated is also reported. The respective results are summarized in Table \ref{hyperreduction}. The $\mathbb{R}\mathbb{E}_{\mathbf{u}}$ error measure of Equation \eqref{errornorm} is depicted for comparison purposes.  

\begin{table}[!ht]
	\caption{Computational efficiency considerations and hyper reduction performance of the pROM. Results are depicted for the small subdomains scenario presented in detail in Table \ref{NarrowRectangle}. All models are simulated with the same number of cores and using a fixed time-step.}
	\label{hyperreduction}
		\begin{tabular} { P{4.5cm} P{2.00cm} P{2.50cm} P{2.75cm} P{2cm} }
		\toprule
		\multicolumn{5}{c}{Hyper reduced pROM}\\ 
		\\
		 & CPU timing \quad (1 core - sec) & Speed-up factor\ \quad(1 core) & Average \break $\mathbb{R}\mathbb{E}_{\mathbf{u}}$ error & \small{Yielding mesh size} \\
		\hline
	High Fidelity Model & $1.19 \times 10^4$& $1.0$& - & 888\\
	pROM approximation & $1.02 \times 10^4$& $1.16$& $\leq1\%$ & 888\\
	Hyper pROM ($\tau=0.01$)  & $3.06 \times 10^2$& 38.80 & 0.0284 & 38\\
	Hyper pROM ($\tau=0.001$)  & $6.02 \times 10^2$& $19.72$& 0.0234 & 75\\
	\bottomrule
	\end{tabular}
\end{table}

\begin{figure*}[!ht]
        \centering
            \centering
            \includegraphics[scale=0.90]{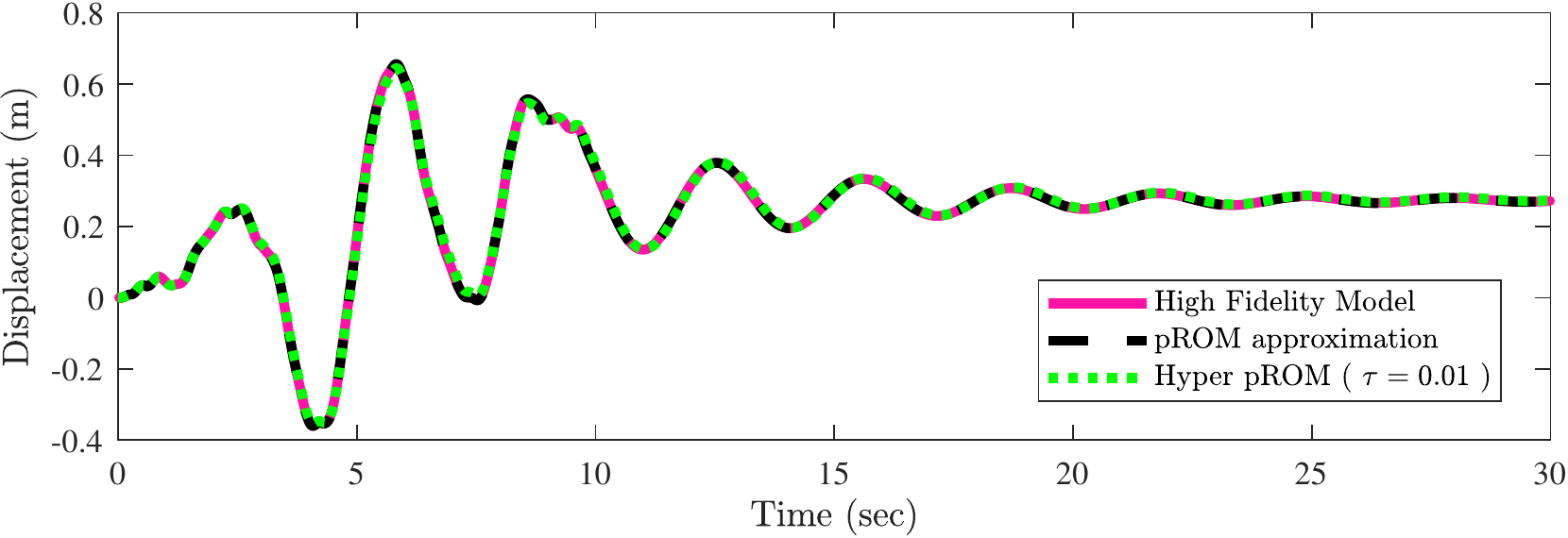}
            \caption[]%
            {Displacement time history of the tip for the hyper pROMs implemented for the example reference region of Table \ref{NarrowRectangle}. The parametric configuration [1.25,36] is presented as validation.}    
            \label{fig:TH}
\end{figure*}

To evaluate both the computational efficiency and the sensitivity of the performance of the hyper-reduction method (ECSW) with respect to the tolerance parameter $\tau$ in Equation \eqref{xi}, Table \ref{hyperreduction} compares four different approaches. Namely the HFM, the pROM approach of this paper without any hyper reduction strategy and two hyper reduced pROMs with a different tolerance parameter. In addition, Figure \ref{fig:TH} presents a comparison on reproducing the time history response on the tip of the tower for an example validation input of [1.25,36].

These results suggest that the hyper reduced pROMs of Table \ref{hyperreduction} perform well in sufficiently capturing the response of the wind turbine tower under parametric excitation input and achieve a considerable computational speed up. This way, the resulting hyper pROM addresses the problem of the accurate modeling of ‘as-is’ nonlinear dynamic structural systems (or components) as defined in section \ref{PrSt}. A sufficient reduced order approximation under variability of the structural properties or of the acting loads is achieved in parallel with significant reduction of the computational load. 

\section{Limitations and Conclusions} \label{5}
This paper implements a physics-based pROM for modeling the dynamic behavior of nonlinear structural systems across a range of parameters, which pertain to i) the structural configuration (material and hysteresis properties), and ii) the spectral and temporal characteristics of the acting loads. A variant to the pROM local bases interpolation technique is here introduced in the context of material nonlinearity. The proposed technique was validated firstly on a simplified shear frame example, where variability of structural parameters is investigated, and secondly on the more complex problem of a wind turbine tower subjected to earthquake loading. The efficiency of the pMOR strategy employed was demonstrated both in terms of displacements and stress states, as stresses are the fundamental quantities of interest in structural monitoring applications. The implemented pROM delivered sufficiently accurate results along with a considerable computational speed up, demonstrating the capabilities of the approach. 

However, certain limitations ought to be acknowledged. Firstly, the implemented pROM was validated on parametric analyses involving a maximum of two parameters. This implies that the interpolation weighting factors, the structure of the grid and the overall scheme are subject to this limitation. In a multi parametric case study these aspects need to be treated with more sophisticated methods. In addition, the interpolation weighting factors are computed based on the distance of the snapshots on the parametric domain. In a more general framework, snapshots have to be clustered based on the premise that their solutions are spanned by the same or similar subspaces. These limitations form challenges to be addressed in future work, along with a more generalized parameterization of earthquake time histories, following the work of \citet{paramearth}.

\section*{Acknowledgements}

\small{This research has been funded from the Sandia National Laboratories, the European Research Council, under the ERC Starting Grant WINDMIL (ERC-2015-StG \#679843) on the topic of  "Smart Monitoring, Inspection and Life-Cycle Assessment of Wind Turbines", as well as from the European Union's Horizon 2020 research and innovation programme under the Marie Sk\l{}odowska-Curie grant agreement No. 795917 ``SiMAero, Simulation-Driven and On-line Condition Monitoring with Applications to Aerospace''.

Sandia National Laboratories is a multimission laboratory managed and operated by National Technology and Engineering Solutions of Sandia, LLC, a wholly owned subsidiary of Honeywell International Inc., for the U.S. Department of Energy’s National Nuclear Security Administration under contract DE-NA0003525.  This paper describes objective technical results and analysis. Any subjective views or opinions that might be expressed in the paper do not necessarily represent the views of the U.S. Department of Energy or the United States Government.}

\bibliographystyle{unsrtnat}
\bibliography{Literature}

\end{document}